\documentclass[preprint,12pt]{elsarticle}

\usepackage{amsmath,amssymb,amsthm}
\usepackage{ifthen,verbatim,here}

\usepackage{mathrsfs}
\usepackage{yhmath}
\usepackage{booktabs}
\usepackage{underscore}
\usepackage[framed,numbered,autolinebreaks,useliterate]{mcode}
\usepackage{longtable}
\usepackage{url}
\usepackage{pdfpages}
\usepackage[english]{babel}
\usepackage{verbatim,here}
\usepackage[T1]{fontenc}
\usepackage{floatflt,graphicx,graphics}
\usepackage{a4wide}

\usepackage{subfig}
\graphicspath{{Figures/}}

\numberwithin{equation}{section}
\nonstopmode
\theoremstyle{plain}

\newtheorem{theorem}[equation]{Theorem}

\theoremstyle{definition}

\newtheorem{defn}[equation]{Definition}

{\qed\bigskip}

\DeclareMathOperator{\capa}{\mathrm{cap}}

\newcounter{alphabet}

\newcommand{\be}{\begin{eqnarray}}
\newcommand{\ee}{\end{eqnarray}}
\newcommand{\ba}{\begin{array}}
\newcommand{\ea}{\end{array}}
\newcommand{\ben}{\begin{eqnarray*}}
\newcommand{\een}{\end{eqnarray*}}

\newcommand{\B}{\mathbb{B}}

\renewcommand{\th}{\,\textnormal{th}}

\newcommand{{\tth}}{\mathrm{th}}
\newcommand{{\sh}}{\mathrm{sh}}
\newcommand{{\ch}}{\mathrm{ch}}

\newcommand{\K}{\mathcal{K}}

\renewcommand{\Im}{{\,\operatorname{Im}\,}}
\renewcommand{\Re}{{\,\operatorname{Re}\,}}

\newcommand {\M} {\mathsf{M}}

\renewcommand{\i}{\mathrm{i}}

\newcommand{\bI}{{\bf I}}
\newcommand{\bM}{{\bf M}}
\newcommand{\bN}{{\bf N}}
\renewcommand{\Im}{{ \rm Im}\,}
\renewcommand{\Re}{{ \rm Re}\,}

\newcommand{\DD}{\mathbb{D}}
\newcommand{\cT}{\mathcal{T}}
\newcommand{\eE}{\varepsilon}
\newcommand{\osc}{\mathrm{osc}}
\font\fFt=eusm10 %scaled 1200
\font\fFa=eusm7  %scaled 1200
\font\fFp=eusm5  %scaled 1200
\def\K{\mathchoice
{\hbox{\,\fFt K}}
{\hbox{\,\fFt K}}
{\hbox{\,\fFa K}}
{\hbox{\,\fFp K}}}
\font\fFt=eusm10 %scaled 1200
\font\fFa=eusm7  %scaled 1200
\font\fFp=eusm5  %scaled 1200

\newcounter{minutes}
\divide\time by 60
\newcounter{hours}
\multiply\time by 60
\addtocounter{minutes}{-\time}
\newtheorem{subsectionB}[equation]{}
\newtheorem{subsubsectionB}[equation]{}

\journal{}

\begin{document}

\begin{frontmatter}
\title{Conformal capacity and
polycircular domains}
%{Conformal capacity and
%multiply connected domains with circular arc polygons as boundaries 
%}
%\title[Conformal capacity and
%multiply connected domains]
%{Conformal capacity and
%multiply connected domains with circular arc polygons as boundaries 
%}

%\def\thefootnote{}
%\footnotetext{
%\texttt{\tiny File:~\jobname .tex,
%          printed: \number\year-\number\month-\number\day,
%          \thehours.\ifnum\theminutes<10{0}\fi\theminutes}
%}
%\makeatletter\def\thefootnote{\@arabic\c@footnote}\makeatother
\author[H. Hakula]{Harri Hakula}
\address[H. Hakula]{Aalto University, Department of Mathematics and Systems Analysis, P.O. Box 11100, FI-00076 Aalto, Finland}
\ead{harri.hakula@aalto.fi}
\author[M.M.S. Nasser]{Mohamed M. S. Nasser}
\address[M.M.S. Nasser]{Mathematics Program, Department of Mathematics, Statistics and Physics, College of Arts and Sciences, Qatar University, Doha, Qatar}
\ead{mms.nasser@qu.edu.qa}
%\author[O. Rainio]{Oona Rainio}
%\address{Department of Mathematics and Statistics, University of Turku, FI-20014 Turku, Finland}
%\email{ormrai@utu.fi}
\author[M. Vuorinen]{Matti Vuorinen}
\address[M. Vuorinen]{Department of Mathematics and Statistics, University of Turku, FI-20014 Turku, Finland}
\ead{vuorinen@utu.fi}

%\subjclass[2010]{Primary !!!; Secondary !!!}
\begin{abstract}

We study numerical conformal mapping of multiply connected planar domains with boundaries consisting 
of unions of finitely many circular arcs, so called polycircular domains. 
We compute the conformal  capacities of condensers defined by polycircular domains. 
Experimental error estimates are provided for the computed capacity and, when possible, 
the rate of convergence under refinement of discretisation is analysed. The main ingredients 
of the computation are two computational methods, on one hand the boundary integral equation 
method combined with the fast multipole method   and on the other hand the $hp$-FEM method. 
The results obtained with these two methods agree with high accuracy.
\end{abstract}
\end{frontmatter}

%%%FILE: sec111.tex
%sec107b.tex
%%%%%%%%%%%%%%%%%%%%%%%%%%%%%%%%%%%%%%%%%
%%%%%%%%%%%%%%%%%%%%%%%%%%%%%%%%%%%%%%%%%
\section{Introduction}
%%%%%%%%%%%%%%%%%%%%%%%%%%%%%%%%%%%%%%%%%
%%%%%%%%%%%%%%%%%%%%%%%%%%%%%%%%%%%%%%%%%
%%%%%%%%%%%%%%%%%%%%%%%%%%%%%%%%%%%%%%%%%

During the past century, \emph{domain functionals} of planar domains have been extensively studied in geometric function theory
and other fields of mathematical analysis. Some of these functionals are harmonic measure, hyperbolic distance and its various
generalizations, conformal radius, and various domain characteristics, for instance, uniform perfectness. One way to classify these notions is
their invariance properties such as invariance under conformal mappings, M\"obius transformations or stretchings. We study here
one of these functionals, the {conformal capacity} of a condenser. A \emph{condenser} is a pair $(G,E)$ where $ G \subset
{\mathbb{R}^2}$ is a domain and $E \subset G$ is a non-empty compact set and the  
\emph{conformal capacity} is defined
as follows
\begin{align}\label{def_condensercap}
{\rm cap}(G,E)=\inf_{u\in A}\int_{G}|\nabla u|^2 \, dm,
\end{align}
where $A$ is the family of all harmonic functions $u$ defined in $G$ with $u(x)\geq1$ for all $x\in E$ and $u(x)\to 0$ when $x\to\partial G$. 
In concrete applications the sets $E$ and $\partial G$ have a simple geometry, both have a finite number of components, each component being a piecewise smooth curve. 
In this paper we assume that $G$ is a simply connected domain and the domain $G\setminus E$ is multiply connected and all of its boundary components are
finite unions of circular arcs or linear segments. Such  domains are known as \emph{polycircular domains}. This term is due to D.~Crowdy~\cite{c}.
In this case, it is known that the infimum 
 is attained by a harmonic function \cite[p. 65]{Ah}, \cite{du}.  Moreover, this extremal function is a solution 
  to the Dirichlet problem
\begin{subequations}\label{eq:u-Dir}
\begin{eqnarray}
\label{eq:u-Dir-1}
\Delta u&=& 0, \quad {\rm on}\;\;G\setminus E, \\
\label{eq:u-Dir-2}
u      &=& 0, \quad {\rm on}\;\;\partial G, \\
\label{eq:u-Dir-3}
u      &=& 1, \quad {\rm on}\;\;\partial E,
\end{eqnarray}
\end{subequations}
and hence the capacity can be expressed in terms of this extremal function as 
% COLOR BEGIN {\color{black}
\begin{equation}\label{eq:cap}
\capa(G,E)=\iint\limits_{\Omega}|\nabla u|^2\, dm
\end{equation}
where $\Omega=G\setminus E$.
%In this paper, we assume $\Omega$ is a multiply connected domain such
%that each of its boundary components consists of a finite number of
%segments and/or circular arcs. Such  domains are known as
%\emph{polycircular domains}. This term is due to Crowdy~\cite{c}. 
Let
$\Gamma_0=\partial G$, $\partial E=\Gamma_1\cup\cdots\cup\Gamma_m$,
where $\Gamma_0$ encloses all the other curves
$\Gamma_1,\ldots,\Gamma_m$. Then each of the boundary components
$\Gamma_j$ is a piecewise smooth Jordan curve with a finite number of
corner points. We assume that none of these corner points is a cusp, that is, the boundary arcs meet at nonzero angles. 
The orientation of the total boundary
$\Gamma=\partial\Omega=\Gamma_0\cup\Gamma_1\cup\cdots\cup\Gamma_m$ is
assumed to be such that $\Omega$ is on the left of $\Gamma$. 
The above Dirichlet problem \eqref{eq:u-Dir} in multiply connected domains is a wide area of research, see e.g. \cite{ais}.

In their classical book  G. Polya and G. Szeg\"o~\cite{ps}
studied extensively various estimates for capacities in terms of various 
domain functionals, however
their capacity was not the conformal capacity. V. Dubinin \cite{du} 
studied conformal capacity in
many function theoretic applications and the book of J. Garnett and 
D. Marshall  \cite{garmar}  is an extensive treatise of harmonic measure 
and its applications to potential theory and geometric function theory.

\begin{figure}
    \centering
    \subfloat[]{\label{fig:AppetiserA}\includegraphics[width=0.4\textwidth]{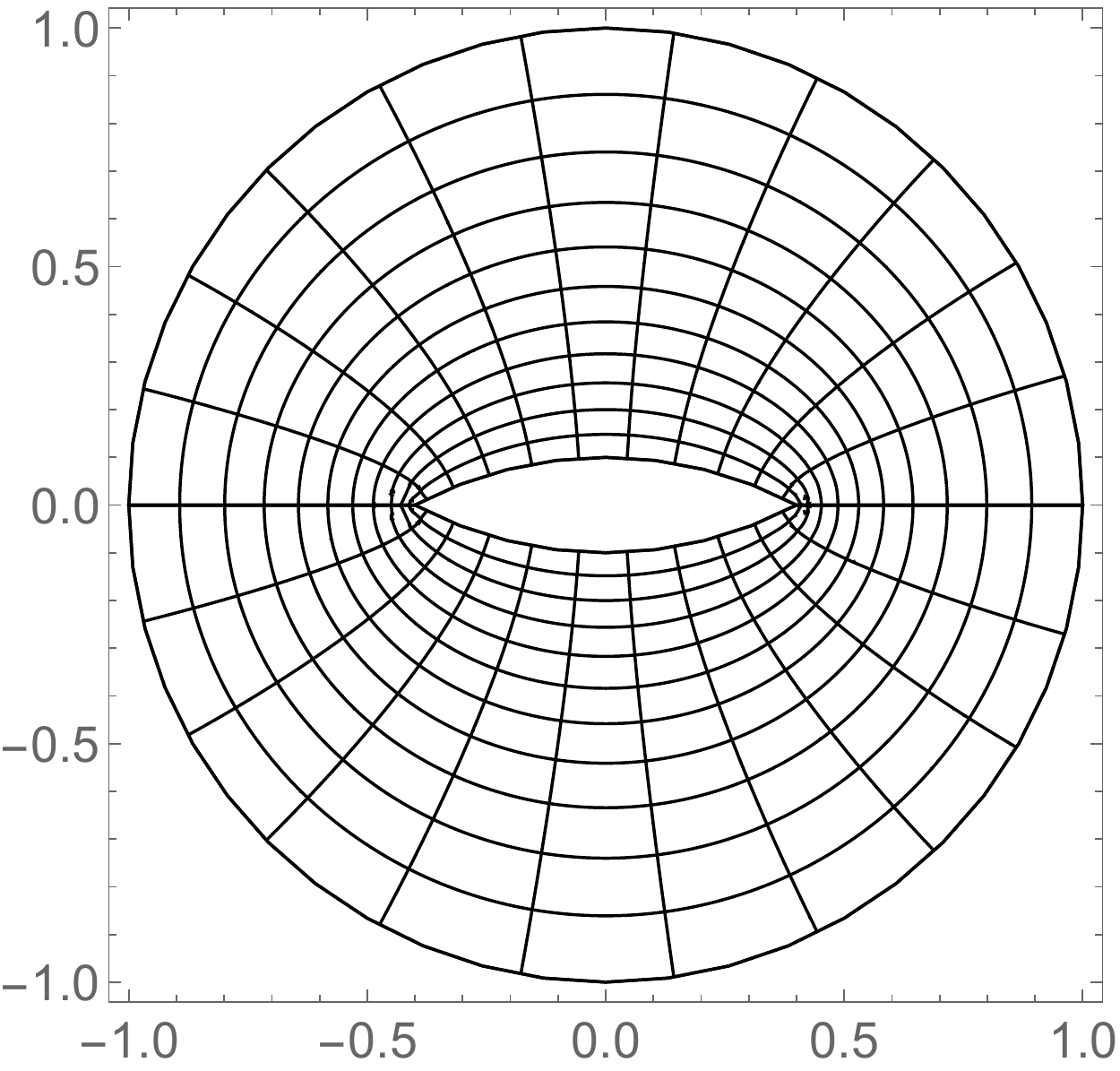}}\quad
    \subfloat[]{\label{fig:AppetiserB}\includegraphics[width=0.45\textwidth]{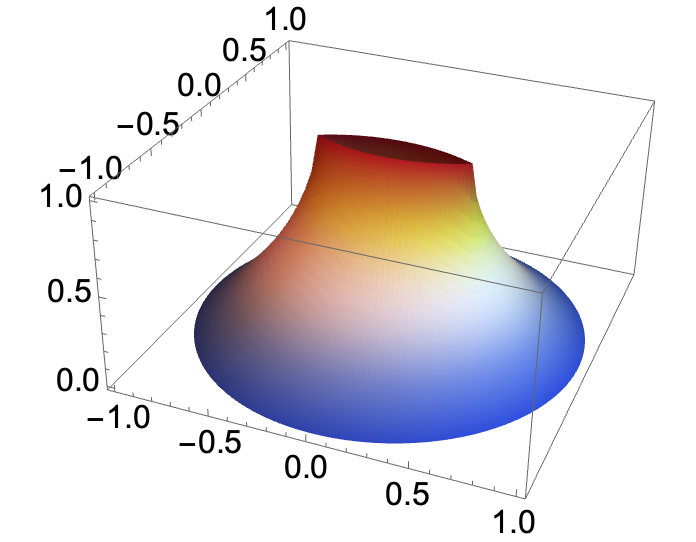}}
 \caption{A polycircular ring domain with the lens boundary $\partial E$ is defined by two circular arcs passing through the points $(-r, 0), (0,\pm s), (r, 0)$, respectively. Here $r= 2/5,\, s=1/10$.
The capacity of this ring domain can be computed by computing a conformal map from this domain onto an annular domain. The left figure  (a) shows the gridlines mapped onto a polar grid, consisting of radial segments
and circles, under this map. The numerical value of the
    capacity given by the $hp$-FEM is ${\rm cap}(\B^2,E) = 4.371029672008615$. In (b)
    the potential function is illustrated.}\label{fig:Appetiser}
    %\caption{Example of a polycircular ring domain. The lens boundary $\partial E$ is defined by two circular arcs passing
    %through the points $(-r, 0), (0,\pm s), (r, 0)$, respectively. Here $r= 2/5,\, s=1/10$. The computed
    %capacity is ${\rm cap}(\B^2,E) = 4.371029672008615$. In (a) the conformal map is shown, and in (b)
    %the potential function is illustrated.}\label{fig:Appetiser}
\end{figure}

In spite of their important role in geometric function theory, explicit formulas or numerical values of conformal invariants are known only in the simplest cases. Roughly 
speaking one can say that when the domain connectivity increases, the more difficult it is to find explicit formulas even for simple domains like the circular domains where all the boundary components are circles. 
The connectivity of $G \setminus E$ 
depends on the number of topological components of $E$ and $G \setminus E$ is doubly connected if $\partial E$ is a simple curve. Our goal here is to continue our earlier work \cite{nrrvwyz,nv1,nv} and to develop algorithms for the computation of ${\rm cap}(G,E)$. For an example, see Figure~\ref{fig:Appetiser} where $G=\B^2$ and $\B^2$ is the unit disk. 

Numerical computation of conformal capacity has been studied in some earlier papers  \cite{nv1,nv}, but we have not seen any results for conformal capacities of condensers for multiply connected polycircular domains of the type described above.
Due to the conformal invariance of the conformal capacity, an  auxiliary step in the computation is often to apply conformal mappings onto a canonical domain to simplify the geometry \cite{kuh}. For example, the conformal capacity for the doubly connected domain shown in Figure~\ref{fig:Appetiser} can be computed by finding the conformal mapping $w=f(z)$ from this domain onto an annulus domain $q<|w|<1$ and hence the conformal capacity is $2\pi/\log(1/q)$. 
The books of N. Papamichael- N. Stylianopoulos \cite{ps10}  and P. Kythe \cite{ky} and the long
survey of R. Wegmann \cite{Weg05} are valuable general sources, see also \cite[pp.14-16]{ps10} for
an overview of the literature.
We will now review the most recent literature from the point of view of polycircular domains.
 
Conformal mapping of Jordan domains with the boundary consisting of a union of finitely
 many circular arcs, has been studied by P. Brown and M. Porter \cite{bp,p}, 
 by U. Bauer and W. Lauf \cite{BL}, and, in particular, by D. Crowdy~\cite{c}, D. Crowdy and A. Fokas~\cite{c2}, and by D. Crowdy, A. Fokas, and C. Green~\cite{c2}.  See also \cite{BjG,Ho,Tr}.
In~\cite{c,c2,c3}, D. Crowdy and his coauthors developed an extensive theory based on
 Schottky's prime function method for the construction of conformal mappings from multiply connected circular domains onto polycircular arc domains. 
In the case of multiply connected polycircular domains, conformal mappings onto canonical domains have been investigated also in \cite{ber,gol,koe,neh,tsu}.
Several types of canonical domains are used as an intermediate step of computation.
Some examples are those where the boundary components are circles or parallel slits or concentric circular arcs. M. Badreddine, T. DeLillo, and S. Sahraei~\cite{bds} compare  several numerical conformal mapping methods for multiply connected polycircular domains.
Computer graphics applications of conformal mapping of multiply connected domains appear in~\cite{kyyg}.
The goal of our research here is different from the
above references as our focus is on conformal capacity. For this purpose we apply two
numerical methods, the $hp$-FEM and the method based on the Boundary Integral Equation (BIE) with the generalized Neumann kernel.
Using these methods, our primary goal is to present detailed
numerical reference results on capacities of multiply connected polycircular domains that have not been considered before. The secondary goal is to verify the correctness of the BIE method by taking the hp-FEM, a well-established standard method, as the reference and to compare the performance of these two methods.

%We give new numerical methods for this same case and present our results in the form %of numerical tables, graphics, and analysis of algorithm performance.

The quantity of interest, the capacity of a condenser, is directly applicable in
engineering contexts. 
Since the Dirichlet problem is one of the primary numerical model problems, any standard solution technique can be viewed as having been validated. In this paper, the novelty is the method of verification of the results. Instead of merely comparing different
discretizations, we make use of \textit{two different methods} to verify the correctness of the results.
The first method is based on the BIE with the generalized
 Neumann kernel as developed and implemented in MATLAB in a series of papers 
 during the past two decades, see e.g. \cite{Nas-Siam1}-\cite{Nas-PlgCir}. 
 The method uses the fast multipole method implementation from \cite{Gre-Gim12} 
 for the speed-up of solving linear equations with a special structure. The BIE method has been used to solve several problems 
in domains of very high connectivity, domains with piecewise smooth boundaries, 
domains with close-to-touching boundaries, and in domains of real-world problems. In a recent 
 series of papers (\cite{nrv2}-\cite{nv}), the method was applied for the 
 capacity computation of planar condensers and for the analysis of isoperimetric problems for capacity. In particular, we will make use of the very recent results~\cite{nv1}.

  The second method, $hp$-FEM, is based on a
Mathematica implementation developed and widely tested by H. Hakula during the past
two decades, \cite{HRV,HRV2}. The method allows one to incorporate \textit{a priori}
knowledge of the singularities into highly non-quasiuniform meshes, 
within which the polynomial order can vary from element to element. In the class of problems
discussed here, one would expect exponential convergence in the natural norm if 
the discretisation is refined properly.

%%%% BEGIN COLOR {\color{red}
 
The BIE method converges exponentially for analytic boundaries and algebraically 
for piecewise smooth boundaries, see e.g.,~\cite{Nas-ETNA}.
However, so far, no numerical comparison with other methods has been made in the literature, especially for domains with corners. 
As mentioned above, one of the objectives of this paper is to present such a comparison where the BIE is compared against the $hp$-FEM.
We obtain surprisingly good agreements between the two methods. The 
obtained results show that the BIE indeed gives accurate results for domains 
with corners even if the angles at the corners are small. 
%%%% END COLOR }

The structure of this paper is as follows. In Section 2, we introduce some notation
and terminology. Section 3 is a short summary of the two methods we use, the $hp$-FEM
and the boundary integral equation method. In Section 4 we compute the capacities of
several multiply connected polycircular condensers and compare the performance of
the methods. Moreover, the computational error is analyzed as a function of the
number of degrees of freedom. Section 5 deals with condensers $(\B^2,E)$ where the
compact set $E$ has lens-shaped structure and we study how closely the theoretical
M\"obius invariance can be observed in numerical computations. Some concluding remarks are given in Section 6.

%FILE: sec206.tex
%%%%%%%%%%%%%%%%%%%%%%%%%%%%%%%%%%%%%%%%%
%%%%%%%%%%%%%%%%%%%%%%%%%%%%%%%%%%%%%%%%%
\section{Preliminary notions}
%%%%%%%%%%%%%%%%%%%%%%%%%%%%%%%%%%%%%%%%%
%%%%%%%%%%%%%%%%%%%%%%%%%%%%%%%%%%%%%%%%%
%%%%%%%%%%%%%%%%%%%%%%%%%%%%%%%%%%%%%%%%%

The capacity of a condenser defined in the introduction \eqref{def_condensercap} can be
defined in terms of the Dirichlet problem \eqref{eq:u-Dir} as well as in many other 
equivalent ways as shown in \cite{du}, \cite{GMP}. 
First, the family $A$ of functions in \eqref{def_condensercap}   
may be replaced by several other families by \cite[Lemma 5.21, p. 161]{GMP}. Furthermore, the capacity  is
equal to the modulus of a curve family 
\begin{equation} \label{capmod}
\capa(G,E)=\M(\Delta(E,\partial G;G)),    
\end{equation}
where $\Delta(E,\partial G;G)$ is the family of all curves joining $E$ with the boundary 
$\partial G$ in the domain $G$ and $\M$ stands for the modulus of a curve family 
\cite[Thm 5.23, p. 164]{GMP}. For the basic facts about capacities and moduli, 
the reader is referred to \cite{du, GMP, HKV}.

\subsectionB{\textit{Quadrilaterals}\\}
%\subsection{Quadrilaterals}
A Jordan domain in the complex plane $\mathbb C$ is a domain with boundary homeomorphic to the
unit circle. A {\it quadrilateral} is a Jordan domain $D$ together with four distinguished points $z_1,z_2,z_3,z_4 \in \partial D$ which define a positive orientation of the boundary. In other words, if we traverse the boundary,  then the points occur in the order of indices and the domain $D$ is on the left-hand side. The quadrilateral
is denoted by $(D;z_1,z_2,z_3,z_4).$ The {\it modulus of the quadrilateral}  \cite[ p. 52]{ps10}
is a unique positive number $h$ such that $D$ can be conformally mapped by some conformal map $f$ onto the rectangle with vertices
$0,1,1+\i h, \i h$ with
$$f(z_1) = 0, \quad f(z_2) =1,  \quad f(z_3) =1+\i h ,  \quad f(z_1) =\i h\,.$$
The modulus is denoted ${\rm mod}(D;z_1,z_2,z_3,z_4).$  The following basic formula is often used:
\begin{equation}\label{reciprel}
{\rm mod}( D;z_1,z_2,z_3,z_4) \,{\rm mod}( D;z_2,z_3,z_4, z_1)=1\,.
\end{equation}
A rectangle with sides $a$ and $b$ and its vertices define a quadrilateral. 
Depending on the labeling of its vertices, its modulus is either $a/b$ or $b/a.$ 
An alternative equivalent definition of the modulus of a quadrilateral is based on the mixed Dirichlet-Neumann problem by L.V. Ahlfors \cite[Thm 4.5, p. 65]{Ah}.
%\end{nonsec}

\subsectionB{\textit{Quadrilateral modulus and curve families}\\}
%\subsection{Quadrilateral modulus and curve families}
The modulus of a quadrilateral $(D;z_1,z_2,z_3,z_4)$  is connected
with the modulus of the family of all curves in $D,$ joining the opposite boundary arcs $(z_2,z_3)$ and $(z_4,z_1),$ in a very simple way, as follows
\begin{equation} \label{2moduli}
{\rm mod}(D;z_1,z_2,z_3,z_4) = \M(\Delta((z_2,z_3), (z_4,z_1);D)) \,.
\end{equation}
%\end{nonsec}

Let ${\rm mod_{cal}}( D;z_1,z_2,z_3,z_4)$ be the computed approximate value to ${\rm mod}( D;z_1,z_2,z_3,z_4)$. Then, for
every quadrilateral the so-called reciprocal relation (\ref{reciprel})
can be rephrased as an error estimate
(Definition~\ref{def:recip}) which will be used below as one of our experimental error estimates.
\begin{defn}{(Reciprocal Identity and Error)} \label{def:recip} 
We shall call 
\begin{equation} \label{myRecipErr}
\varepsilon_R = |1-{\rm mod_{cal}}( D;z_1,z_2,z_3,z_4) \,{\rm mod_{cal}}( D;z_2,z_3,z_4, z_1)|\,
\end{equation}
 the error measure and $$\varepsilon_N
=\left|\lceil\log_{10}|\varepsilon_R |\rceil\right|$$ the related
error number. 
This error characteristic was systematically used in \cite{HRV,HRV2} 
where it was shown to be usually compatible with other error characteristics.
\end{defn}

%FILE: sec316.tex
%%%%%%%%%%%%%%%%%%%%%%%%%%%%%%%%%%%%%%%%%
%%%%%%%%%%%%%%%%%%%%%%%%%%%%%%%%%%%%%%%%%
%\section{\color{red}Circular-arc polygons}
\section{Numerical methods used in computing capacities on polycircular domains}

In this section the two methods used in this paper, the $hp$-FEM and the
boundary integral equation (BIE) method, are introduced.
Much of the material is standard, however, terminology and error concepts
used in the discussion on numerical experiments is defined here.

\subsectionB{\textit{High-Order Finite Elements}\\}
%\subsection{High-Order Finite Elements}

The finite element method (FEM) is the standard numerical method
for solving elliptic partial differential equations. Since FEM is
an energy minimization method it is eminently suitable for
problems involving Dirichlet energy. In the context of this paper
where the focus is on domains with singularities at the
vertices, the $hp$-FEM variant is the most efficient one
\cite{bs,schwab}. With proper grading of the meshes
even with uniform polynomial order
exponential convergence can be
achieved even in problems with strong corner singularities.

In this section we give a brief overview of the method and our
implementation \cite{HRV,ht}. Of particular importance is the
possibility to estimate the error in the computed quantity of
interest. 
%For quadrilaterals there exists a natural error
%estimate, the so-called \textit{reciprocal relation} which is a
%necessary but not sufficient condition for convergence. However, if the reciprocal
%relation is coupled with 
%\textit{a posteriori} error estimates, we
%can trust the results with high confidence \cite{hno}.
As noted above, in the absence of exact solutions, within the method itself
we have the option of \textit{a posteriori} error estimation \cite{hno}.

The equality (\ref{eq:cap}) shows that the capacity of a
condenser is the Dirichlet integral, i.e., the $H^1$-seminorm
of the potential $u$ squared, or, the energy norm squared, a
quantity of interest which is natural in the FEM setting.

\subsubsectionB{\textit{Mesh refinement and exponential convergence}\\}
\label{sec:refinement}

In the $h$-version the accuracy of the discretization is 
controlled by the sizes of the elements. The degrees of freedom are associated with the nodes of the mesh and
the nodal shape functions induce a partition of unity.

The idea behind the $p$-version is to associate degrees of freedom
to topological entities of the mesh  and control accuracy via the polynomial order.
The shape functions
are based on suitable orthogonal polynomials and their supports
reflect the related topological entity, nodes, edges, faces (in
3D), and interior of the elements. 
For instance, on a triangle of order $p$, there are three nodal shape functions, $p-1$ shapes on each edge,
and $(p-1)(p-2)/2$ shapes in the interior.
If only nodal shape functions are included, the $p$-version reduces to the classical $h$-version.
The $hp$-version simply refers to combination of the two refinement strategies.

For the Dirichlet problem \eqref{eq:u-Dir} it can be shown that if the mesh is graded
appropriately the method converges exponentially in some
general norm such as the $H^1$-seminorm. Moreover, due to the
construction of shape functions, it is natural to have large
curved elements in the mesh without significant increase in the
discretization error. Since the number of elements can be kept
relatively low given that additional refinement can always be
added via elementwise polynomial degree, variation in the boundary
can be addressed directly at the level of the boundary
representation in some exact parametric form.

\begin{figure}
    \centering
    \subfloat[After sixteen levels of refinement.]{\label{fig:CNGA}\includegraphics[width=0.45\textwidth]{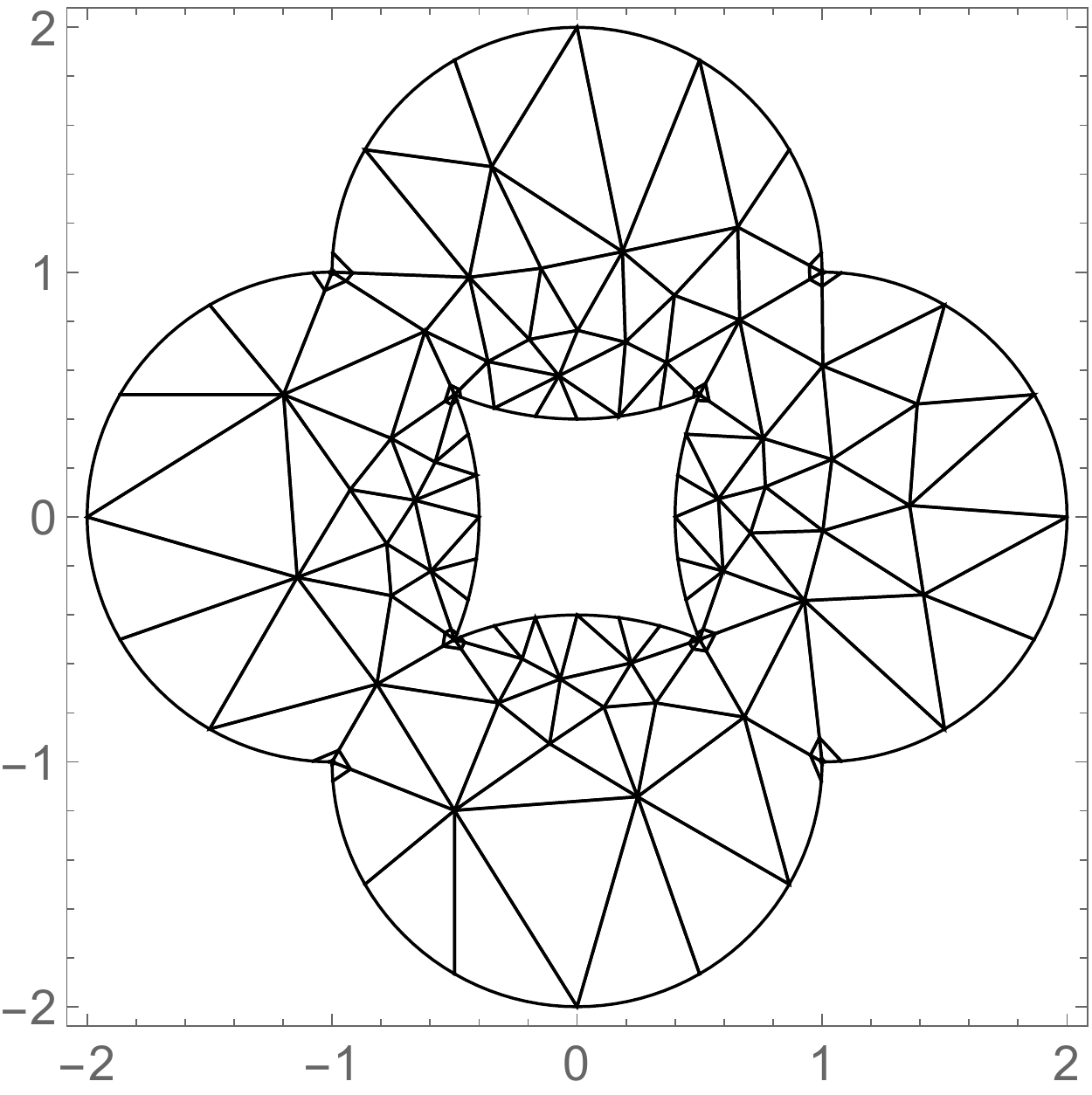}}\quad
    \subfloat[Mesh detail.]{\includegraphics[width=0.45\textwidth]{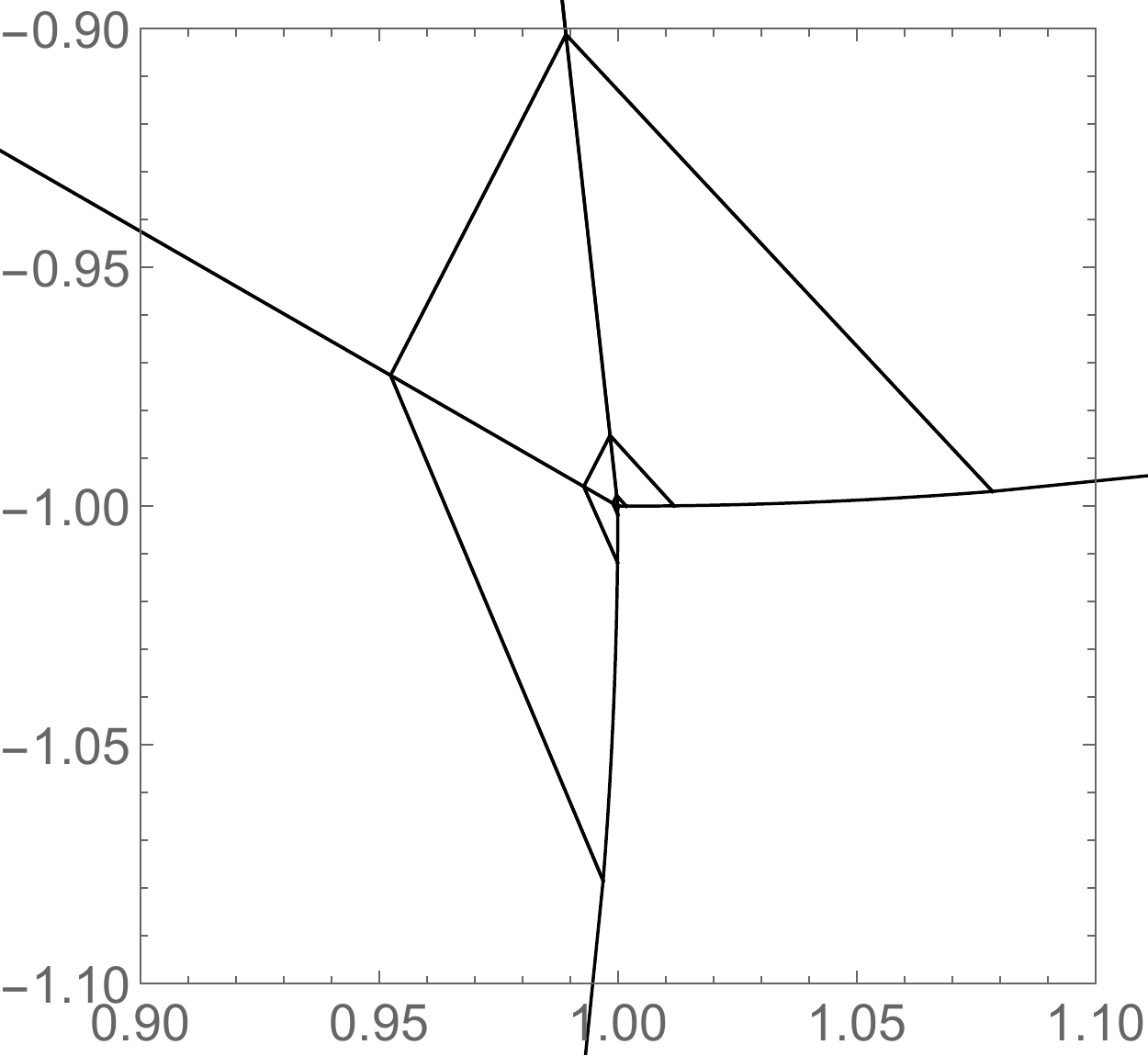}}
    \caption{Example of geometric grading of the mesh. We start with a coarse background triangulation.
    Once the singularities have been identified a priori, a sequence of geometric refinements is
    applied. Notice in (b) that the resulting mesh can include quadrilaterals and triangles.
    The high aspect ratios are compensated with high order polynomial basis.}\label{fig:CNG}
\end{figure}

To fully realize the potential of the $p$-version, one has to
grade the meshes properly and therefore we really use the
$hp$-version here. 
%Consider the meshes in Figures~\ref{fig:CNG}
%and \ref{fig:CQ}. 
In Figure~\ref{fig:CNG} the basic refinement
strategy is illustrated. We start with an initial mesh, where the
corners with singularities are \textit{isolated}, that is, the
subsequent refinements of their neighboring elements do not
interfere with each other. Then the mesh is refined using
successive applications of replacement rules.

In our implementation the geometry can be described in exact
arithmetic and therefore there are not any fixed limits on the
number of refinement levels. In the case of graded meshes one has
to resolve the question of how to set the polynomial degrees at
every element, indeed, a form of refinement of its own. One option
in the case of strong singularities is to set the polynomial
degree based on graph distance from the singularity.
Here, however, the degree $p$ is kept constant over the whole mesh
despite the grading.

The following theorem due to Babu{\v{s}}ka and Guo \cite{bg1}
sets the limit to the rate of convergence. Notice that construction of
appropriate spaces is technical. For rigorous treatment of the theory involved
see Schwab \cite{schwab} and references therein. 
For the Dirichlet problem (\ref{eq:u-Dir})
the following theorem relates the convergence of the capacity (or energy) and the number of degrees
of freedom $N$.
\begin{theorem} \label{propermesh}
Let $\Omega \subset \mathbb{R}^2$ be a polygon, $v$ the FEM-solution of \eqref{eq:u-Dir}, and
let the weak solution $u_0$ be in a suitable countably normed space where
the derivatives of arbitrarily high order are controlled.
Then
\[
\inf_v \|u_0 - v\|_{H^1(\Omega)} \leq C\,\exp(-b \sqrt[3]{N}),
\]
where $C$ and $b$ are independent of $N$, the number of degrees of freedom.
Here $v$ is computed on a proper geometric mesh, where the order of an individual
element is set to be its element graph distance to the nearest singularity.
(The result also holds for meshes with constant polynomial degree.)
\end{theorem}

\subsubsectionB{\textit{Error Estimation: Reciprocal Error and Auxiliary Space Error Estimate}\\}
%\subsubsection{Auxiliary Space Error Estimate}

Assuming that the exact capacity is not known, we have two types
of error estimates available: the reciprocal estimate and an a
posteriori estimate. Naturally, if the exact value is known, we
can measure the true error.

The reciprocal error estimate introduced in Section 2 is rather unusual in the sense
that it is based on \textit{physics}, yet only necessary.
In the context of this paper, the reciprocal error estimator \eqref{myRecipErr} is useful
in cases where the original problem can be reduced to one on a quadrilateral,
for instance, using symmetries.

Our a posteriori error estimator of choice,
the auxiliary space error estimate,
belongs to a family of
hierarchical error estimators. It is eminently suited to be used
in conjunction of $hp$-FEM.
Consider the abstract problem setting with $u$ defined on the standard
piecewise polynomial finite element space on some discretization
$T$ of the computational domain $D$. Assuming that the exact
solution $u \in H_0^1(D)$ has finite energy, we arrive at the
approximation problem: Find $\hat{u} \in V$ such that
\begin{equation}\label{eq:approximation}
  a(\hat{u},v)=l(v)\ (= a(u,v))\quad (\forall v \in V),
\end{equation}where
$a(\,\cdot\,,\,\cdot\,)$ and $l(\,\cdot\,)$, are the bilinear form
and the load potential, respectively. Additional degrees of
freedom can be introduced by enriching the space $V$. This is
accomplished via introduction of an auxiliary subspace or ``error
space'' $W \subset H_0^1(D)$ such that $V \cap W = \{0\}$. We can
then define the error problem: Find $\varepsilon \in W$ such that
\begin{equation}\label{eq:error}
  a(\varepsilon,v)=l(v)- a(\hat{u},v) (= a(u-\hat{u},v))\quad (\forall v \in W).
\end{equation}
This is simply a projection of the residual to the auxiliary space.
In 2D the space $W$, that is, the additional unknowns, can be
associated with element edges and interiors. Thus, as mentioned above, for
$hp$-methods this kind of error estimation is natural. The main
result on this kind of estimators for the Dirichlet problem (\ref{eq:u-Dir}) is the following
theorem.%~\ref{KeyErrorThm}.
\begin{theorem}[\cite{hno}]\label{KeyErrorThm}
There is a constant $K$ depending only on the dimension $d$,
polynomial degree $p$, continuity and coercivity constants $C$ and
$c$, and the shape-regularity of the triangulation $\cT$ such that
\begin{align*}
\frac{c}{C}\,\|\eE\|_1\leq\|u-\hat{u}\|_{1}\leq K
\left(\|\eE\|_{1}+\osc(R,r,\cT)\right),
\end{align*}
where the residual oscillation depends on the volumetric and face
residuals $R$ and $r$, and the triangulation $\cT$.
\end{theorem}

The solution $\varepsilon$ of (\ref{eq:error}) is called the
\textit{error function}. It has many useful properties for both
theoretical and practical considerations. In particular, the error
function can be numerically evaluated and analyzed for any finite
element solution. By construction, the error function is
always zero at the mesh points. In the examples below, the
space $W$ contains edge shape functions of degree $p+1$ and
internal shape functions of degrees $p+1$ and $p+2$.
This choice is not arbitrary but based on careful cost analysis \cite{hno}.

A key measure of the quality of the estimator is its effectivity in a norm of interest,
here we define the effectivity index $\lambda$ as follows:
\begin{equation}
\lambda = \|u-\hat{u}\|_{1} / \|\eE\|_1.
\end{equation}
 
\subsectionB{\textit{The boundary integral equation method}\\}
%\subsection{The boundary integral equation method}
Using Green's formula~\cite[p.~4]{du}, and in view of~\eqref{eq:u-Dir-2}--\eqref{eq:u-Dir-3}, formula~\eqref{eq:cap} can be written as
\begin{equation}\label{eq:cap-n}
\capa(C)=\int_{\Gamma}u\frac{\partial u}{\partial{\bf n}}ds
=\sum_{k=1}^{m}\int_{\Gamma_k}\frac{\partial u}{\partial{\bf n}}ds.
\end{equation}
Further, the harmonic function $u$ is the real part of an analytic function $F$ in $G$ which is not necessarily single-valued. 
Assume for each $k=1,2,\ldots,m$ that $\alpha_k$ is an auxiliary point enclosed by $\Gamma_k$. Then the function $F$ can be written as~\cite{Gak,garmar,Mik64}
\begin{equation}\label{eq:F-u}
F(z)=g(z)-\sum_{k=1}^{m} a_k\log(z-\alpha_k)
\end{equation}
where $g$ is a single-valued analytic function in $G$ and $a_1,\ldots,a_{m}$ are undetermined real constants such that~\cite[\S31]{Mik64}
\begin{equation}\label{eq:ak}
a_k=\frac{1}{2\pi}\int_{\Gamma_k}\frac{\partial u}{\partial{\bf n}}ds, \quad k=1,2,\ldots,m.
\end{equation}
Thus, it follows from~\eqref{eq:cap-n} that
\begin{equation}\label{eq:cap-a}
\capa(C)=2\pi\sum_{k=1}^{m}a_k.
\end{equation}
See~\cite{nv1} for more details. 

The problem~\eqref{eq:u-Dir} above is a particular case of the problem considered 
in~\cite[Eq.~(4)]{nv1}. So, the constants $a_1,\ldots,a_m$ in~\eqref{eq:cap-a} 
can be computed using the BIE method presented in~\cite[Theorem 4]{nv1} which is based on using the BIE with the generalized Neumann kernel. 
For domains with smooth boundaries, the BIE can be solved accurately using the Nystr\"om method with equidistant trapezoidal rule. 
For smooth boundaries of class $C^{q+2}$ and smooth integrand of class $C^q$, the trapezoidal Nystr\"om method converges with order $O(1/n^q)$ where $n$ is the number of mesh points~\cite{kre90}. A MATLAB function {\tt fbie} for solving the BIE with the generalized Neumann using the trapezoidal Nystr\"om method is presented in~\cite{Nas-ETNA} and this function will be used in this paper.
However, in this paper, we consider polycircular domains. For such domains, the solution of the BIE has a singularity in its derivative in the vicinity of the corner points~\cite[p.~390]{atk} and this causes that the equidistant trapezoidal rule yields only poor convergence~\cite{kre90}. 
To achieve a satisfactory accuracy, we discretize the BIE using a graded mesh and then applying the Nystr\"om's method~\cite{atk,kre90,kre91}. To describe such a graded mesh method, we assume that the boundary component $\Gamma_j$ consists of $\ell$ circular arcs or segments, then $\Gamma_j$ has $\ell$ corner points. 
We first parametrize each boundary component $\Gamma_j$ by a $2\pi$-periodic function $\zeta_j(t)$ for $t\in J_j=[0,2\pi]$. The function $\zeta_j(t)$ is assumed to be smooth with $\zeta'_j(t)\ne0$ for all values of $t\in J_j$ such that $\zeta_j(t)$ is not a corner point. 
We assume that $\zeta'_j(t)$ has only the first kind discontinuity at these corner points. If $\zeta_j(\hat t)$ is a corner point, we define $\zeta'_j(\hat t):=\zeta'_j(\hat t-0)$.
Let $J$ be the disjoint union of the $m+1$ intervals $J_j=[0,2\pi]$, $j=0,1,\ldots,m$. We define a parametrization of the whole boundary $\Gamma$ on $J$ by (see~\cite{Nas-ETNA} for the details)
\[
\zeta(t)=\left\{
\begin{array}{cc} 
\zeta_0(t), & t\in J_0, \\ 
\zeta_1(t), & t\in J_1, \\
\vdots \\
\zeta_m(t), & t\in J_m. \\ 
\end{array}
\right.
\]
Then, as noted in the first paragraph in~\cite{kre91}, using the graded mesh method suggested in~\cite{kre90} for discretizing the BIE is equivalent to parameterizing the boundary $\Gamma$ by 
\[
\eta(t)=\zeta(\delta(t)),
\]
where the function $\delta(t)$ is defined in~\cite[pp.~696--697]{LSN17} which is chosen to remove the discontinuity in the derivatives of the solution of the BIE. Then, the BIE can be solved using the MATLAB function {\tt fbie} as in the case of smooth domains.

With the parametrization $\eta(t)$ of the whole boundary $\Gamma$, we define a complex function $A$ by
\begin{equation}\label{eq:A}
A(t) = \eta(t)-\alpha,
\end{equation}
where $\alpha$ is a given point in the domain $\Omega$. The generalized Neumann kernel $N(s,t)$ is 
defined for $(s,t)\in J\times J$ by
\begin{equation}\label{eq:N}
N(s,t) :=
\frac{1}{\pi}\Im\left(\frac{A(s)}{A(t)}\frac{\eta'(t)}{\eta(t)-\eta(s)}\right).
\end{equation}
We define also the following kernel 
\begin{equation}\label{eq:M}
M(s,t) :=
\frac{1}{\pi}\Re\left(\frac{A(s)}{A(t)}\frac{\eta'(t)}{\eta(t)-\eta(s)}\right),
 \quad (s,t)\in J\times J.
\end{equation}
The kernel $N(s,t)$ is continuous and the kernel $M(s,t)$ is singular where the 
singular part involves the cotangent function. Hence, the integral operator $\bN$ 
with the kernel $N(s,t)$ is compact and the integral operator $\bM$ with 
the kernel $M(s,t)$ is singular. Further details can be found in~\cite{Weg-Nas}.

The next theorem follows from~\cite[Theorem 4]{nv1}.

\begin{theorem}\label{thm:method}
For each $k=1,2,\ldots,m$, let the function $\gamma_k$ be defined by
\begin{equation}\label{eq:gam-k}
\gamma_k(t)=\log|\eta(t)-\alpha_k|,
\end{equation}
let $\mu_k$ be the unique solution of the integral equation
\begin{equation}\label{eq:ie}
\mu_k-\bN\mu_k=-\bM\gamma_k,
\end{equation}
and let the piecewise constant function $h_k=(h_{0,k},h_{1,k},\ldots,h_{m,k})$ be given by
\begin{equation}\label{eq:h}
h_k=[\bM\mu_k-(\bI-\bN)\gamma_k]/2.
\end{equation}
Then, the $m+1$ real constants $a_1,\ldots,a_{m},c$ are the unique solution of the linear system
\begin{equation}\label{eq:sys-method}
\left[\begin{array}{ccccc}
h_{0,1}    &h_{0,2}    &\cdots &h_{0,m}      &1       \\
h_{1,1}    &h_{1,2}    &\cdots &h_{1,m}      &1       \\
\vdots     &\vdots     &\ddots &\vdots       &\vdots  \\
h_{m,1}    &h_{m,2}    &\cdots &h_{m,m}      &1       \\
\end{array}\right]
\left[\begin{array}{c}
a_1    \\a_2    \\ \vdots \\ a_{m} \\  c 
\end{array}\right]
= \left[\begin{array}{c}
0 \\  1 \\  \vdots \\ 1  
\end{array}\right].
\end{equation}
\end{theorem}

Note that the function $h_k=(h_{0,k},h_{1,k},\ldots,h_{m,k})$ in~\eqref{eq:h} is a piecewise constant function, i.e., for $k=1,\ldots,m$ and $p=0,1,\ldots,m$, the function $h_k$ is constant on each boundary component $\Gamma_p$ and the real constant $h_{p,k}$ is the value of $h_k$ on $\Gamma_p$.
For each $k=1,2,\ldots,m$, the solution $\mu_k$ of the BIE~\eqref{eq:ie} and the piecewise constant function $h_k$ in~\eqref{eq:h} will be computed using the MATLAB {\tt fbie} from~\cite{Nas-ETNA} by calling
\[
{\tt [muk,hk]=fbie(et,etp,A,gamk,n,iprec,restart,gmrestol,maxit)}
\]
where {\tt et} and {\tt etp} are the discretization vectors of the parametrization $\eta(t)$ of the boundary $\Gamma$ and its derivative. 
Here, $n$ is the number of mesh points on each boundary component $\Gamma_k$, and hence the total number of mesh points on the whole boundary $\Gamma$ is $(m+1)n$. The vectors {\tt et} and {\tt etp} are computed by the MATLAB function {\tt cirarcp3pt} which is based on the above described method for parameterizing the boundary $\Gamma$ and is available in \url{https://github.com/mmsnasser/polycircular}. Then the discretization vectors {\tt A} and {\tt gamk} of the functions $A(t)$ and $\gamma_k(t)$ are computed by~\eqref{eq:A} and~\eqref{eq:gam-k}, respectively.
By computing approximate values of the piecewise constant function $h_k=(h_{0,k},h_{1,k},\ldots,h_{m,k})$ in~\eqref{eq:h}, we obtain the entries of the coefficient matrix of the $(m+1)\times(m+1)$ linear system~\eqref{eq:sys-method}. Since $m+1$ is the number of boundary components of the domain $\Omega=G\setminus E$, we can assume that $m$ is small (say, $m<1000$). Hence, the linear system~\eqref{eq:sys-method} will be solved using the Gauss elimination method. By solving the linear system, we obtain the values of the real constants $a_1,\ldots,a_m$ and hence the capacity $\capa(G,E)$ will be computed by~\eqref{eq:cap-a}.

The MATLAB function {\tt fbie} is based on using the trapezoidal Nystr\"om method. Thus, with the above parametrization $\eta(t)$ of the boundary $\Gamma$, it follows from~\cite[Theorems 2.1 and~4.4]{kre90} that the order of the method is $Cn^{-(2q+1)}$ where $q$ is an integer such that $2q+1\le \hat\alpha p$ where $0<\hat\alpha<1$ and the constant $C$ depends on $\hat\alpha$, $q$, and the integrand. The constant $\hat\alpha$ depends on the integrand which in turns depends on the parametrization $\eta(t)$ (see~\cite[Section~2]{kre90} for more details). Here $p$ is an integer known as the grading parameter and $p\ge2$. It is involved in the definition of the function $\delta(t)$ (see~\cite[pp.~696--697]{LSN17}).
Theoretically, choosing large values of $p$ would increase the order of convergence of the method. However, choosing large value of $p$ means that more grid points will be close to corner points which could cause some of these grids to be numerically equal for large values of $n$ (in double precision calculations). In this paper, we choose $p=3$ which suggests that the accuracy of the method is at best $Cn^{-1}$. However, the numerical results presented below in Figure~\ref{fig:err-BIEN} illustrate that the order of the convergence of the method is better than $Cn^{-1}$.

Finally, a MATLAB function {\tt capm} for computing the capacity using the above described method is given in \url{https://github.com/mmsnasser/polycircular}.

%\end{mysubsection}

%FILE: sec416.tex
%%%%%%%%%%%%%%%%%%%%%%%%%%%%%%%%%%%%%%%%%
%%%%%%%%%%%%%%%%%%%%%%%%%%%%%%%%%%%%%%%%%
%\section{\color{red}Examples} 
\section{Verification experiments on polycircular domains}

\begin{figure}
 \centering
  \subfloat[{$m=1$}]{\label{fig:3-cir-m1}\includegraphics[width=0.3\linewidth]{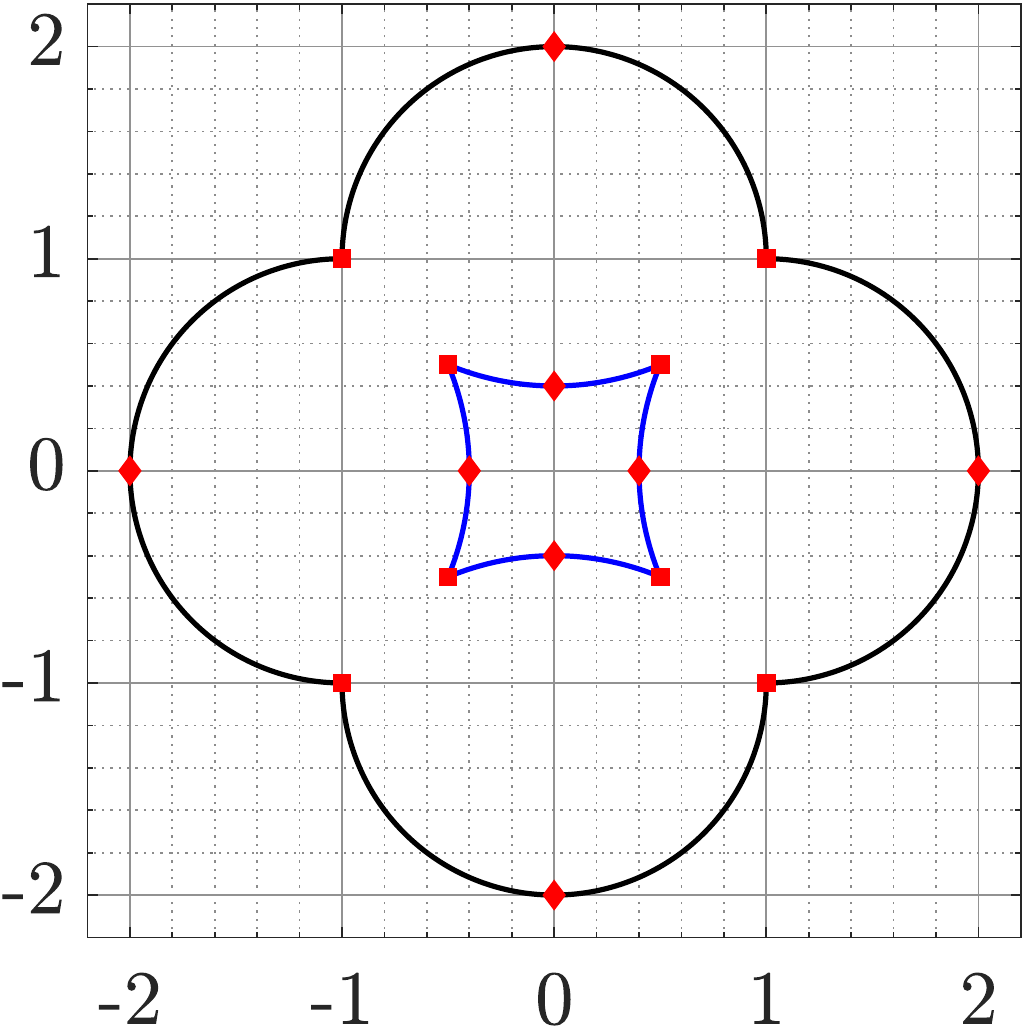}}
  \subfloat[{$m=2$}]{\label{fig:3-cir-m2}\includegraphics[width=0.33\linewidth]{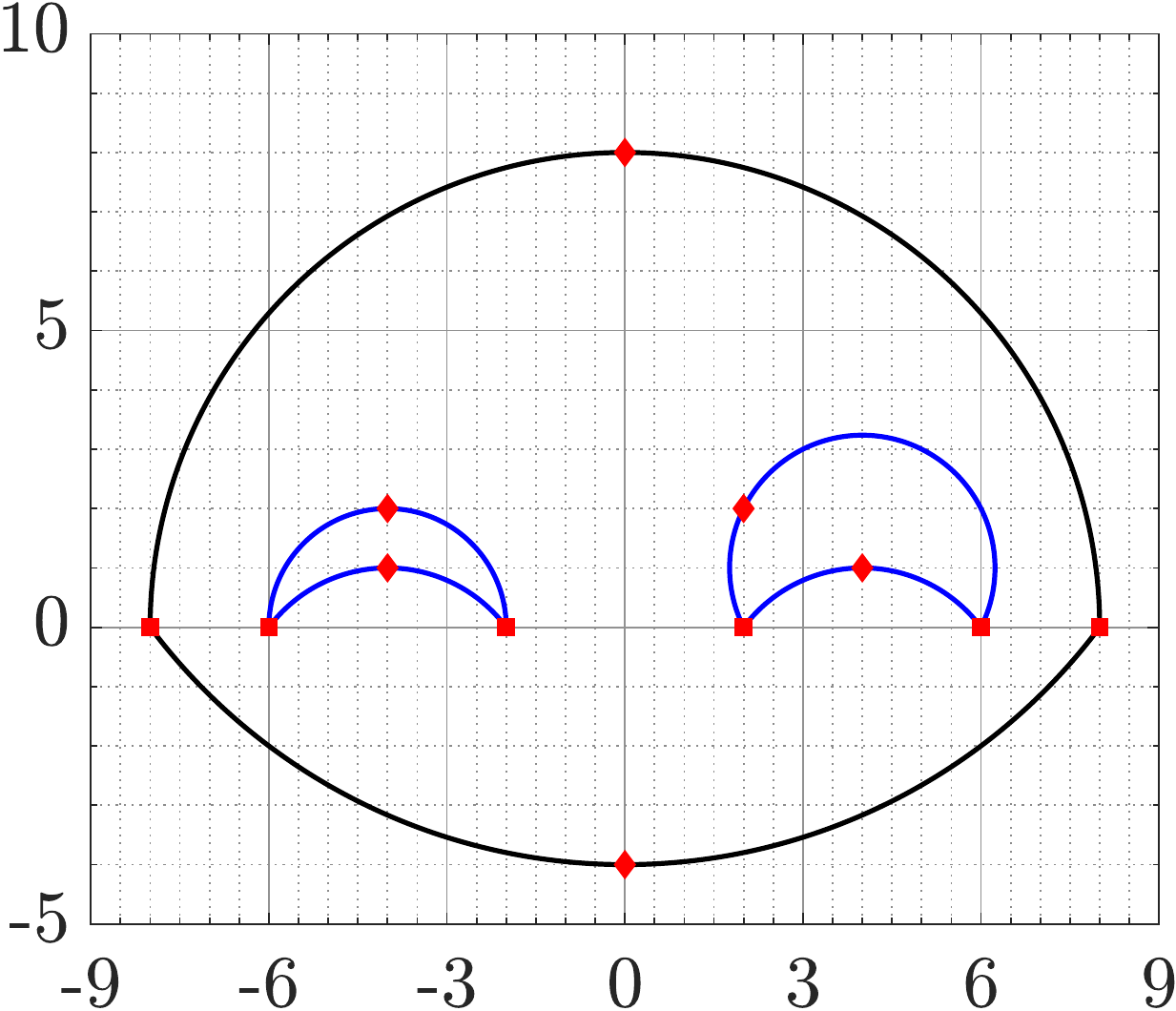}}
  \subfloat[{$m=5$}]{\label{fig:3-cir-m5}\includegraphics[width=0.33\linewidth]{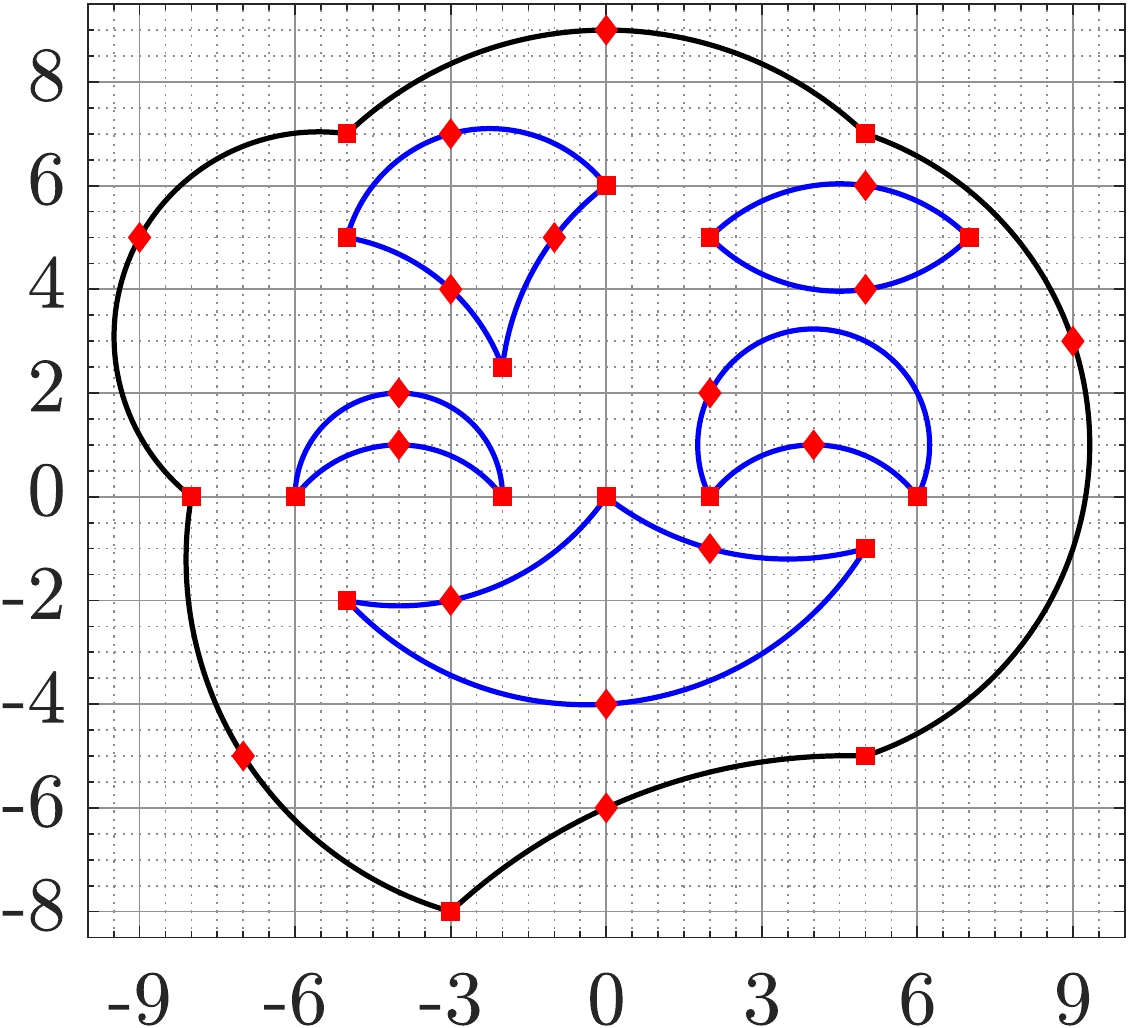}}
  \caption{Multiply connected polycircular domains. 
  The number $m$ refers to the number of interior components.}
  \label{fig:3-cir-mN}
\end{figure}

In this section, we consider three polycircular domains with 
$m=1$, $m=2$, and $m=5$ (see Figure~\ref{fig:3-cir-mN}). 
Each arc is determined by three points, two end points (squares) 
and one interior point (diamond). The purpose of these experiments
is to establish the convergence rates of the methods and their
respective accuracies using  both error estimators and direct comparisons. 
The obtained approximate values of 
$\capa\Omega$ for both methods are shown in Table~\ref{tab:cap-values}.
The reported solution times indicate that 
BIE is generally faster. Whereas the times of BIE depend on the number of mesh points on each boundary component,
on FEM the dependence is on the number of elements and the chosen polynomial order. Moreover, the timings do not reflect
the time spent specifying the geometry, which for FEM can often be the dominant cost measured in time.

%\begin{table}[H]
  %\caption{The values of $\capa\,\Omega$ for the multiply connected polycircular domains in Figure~\ref{fig:3-cir-mN}. For $m=1$ the computed reference value is $5.597545663702171$. Agreement is taken to be the difference between the respective values.}
  %\label{tab:cap-values}
%\centering
%\begin{tabular}{l|l|l|l|l}
    %\hline
        %& BIE       & $hp$-FEM               & Agreement  & Significant digits\\ \hline
%$m=1$   & $5.597545657473350$  & $5.597545663702324$ & $6.23\times10^{-9}$ & 9\\
%$m=2$   & $10.95831566112476$  & $10.95831564100293$ & $6.52\times10^{-9}$ & 10 \\
%$m=5$   & $21.25386110000110$  & $21.25386107560764$ & $2.44\times10^{-8}$ & 9 \\  \hline
%\end{tabular}
%\end{table}

\begin{table}[H]
  \caption{The values of $\capa\,\Omega$ for the multiply connected polycircular domains in Figure~\ref{fig:3-cir-mN}. For $m=1$ the computed reference value is $5.597545663702171$. Agreement is taken to be the difference between the respective values.}
  \label{tab:cap-values}
\centering
\begin{tabular}{|l|c|c||c|c||c|}
    \hline
    & \multicolumn{2}{c||}{BIE Method} & \multicolumn{2}{c||}{$hp$-FEM} & Agreement \\
    \cline{2-5}
	  & Capacity & Time (sec) & Capacity & Time (sec) & \\ \hline
 $m=1$  & $5.597545657473350$ & $2.518$ & $5.597545663702324$ & $19.56$ & $6.23\times10^{-9}$\\
 $m=2$  & $10.95831566112476$ & $5.206$ & $10.95831564100293$ & $12.18$ & $6.52\times10^{-9}$\\
 $m=5$  & $21.25386110000110$ & $25.29$ & $21.25386107560764$ & $69.48$ & $2.44\times10^{-8}$\\
    \hline
\end{tabular}
\end{table}

\begin{figure}
    \centering
    \subfloat[Reciprocal error (logplot). ]{\label{fig:CGA}\includegraphics[width=0.45\textwidth]{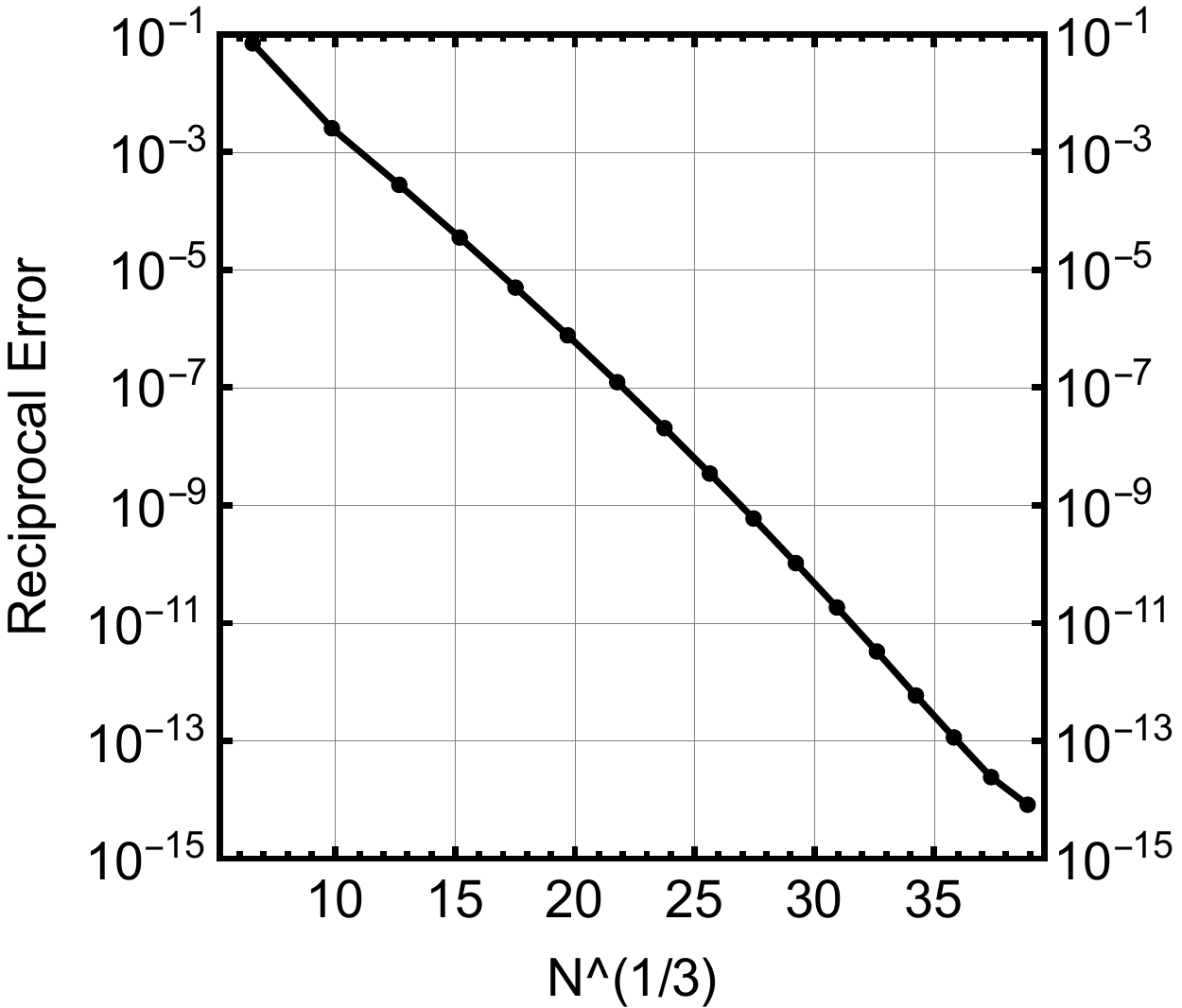}}\quad
    \subfloat[Auxiliary space error (logplot).]{\label{fig:CGB}\includegraphics[width=0.45\textwidth]{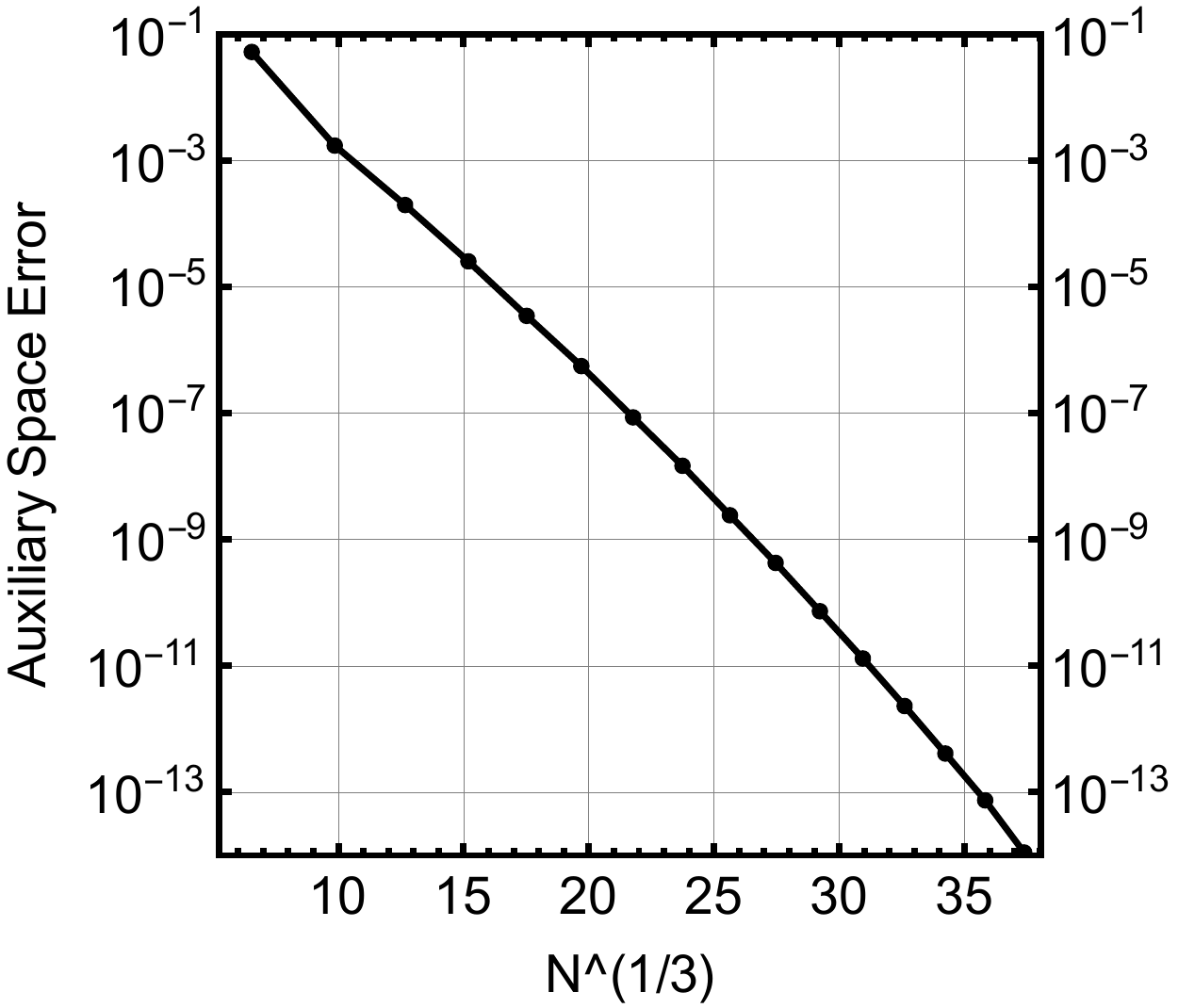}}
    \\
    \subfloat[{Potential $u$.}]{\label{fig:CGC}\includegraphics[width=0.3\textwidth]{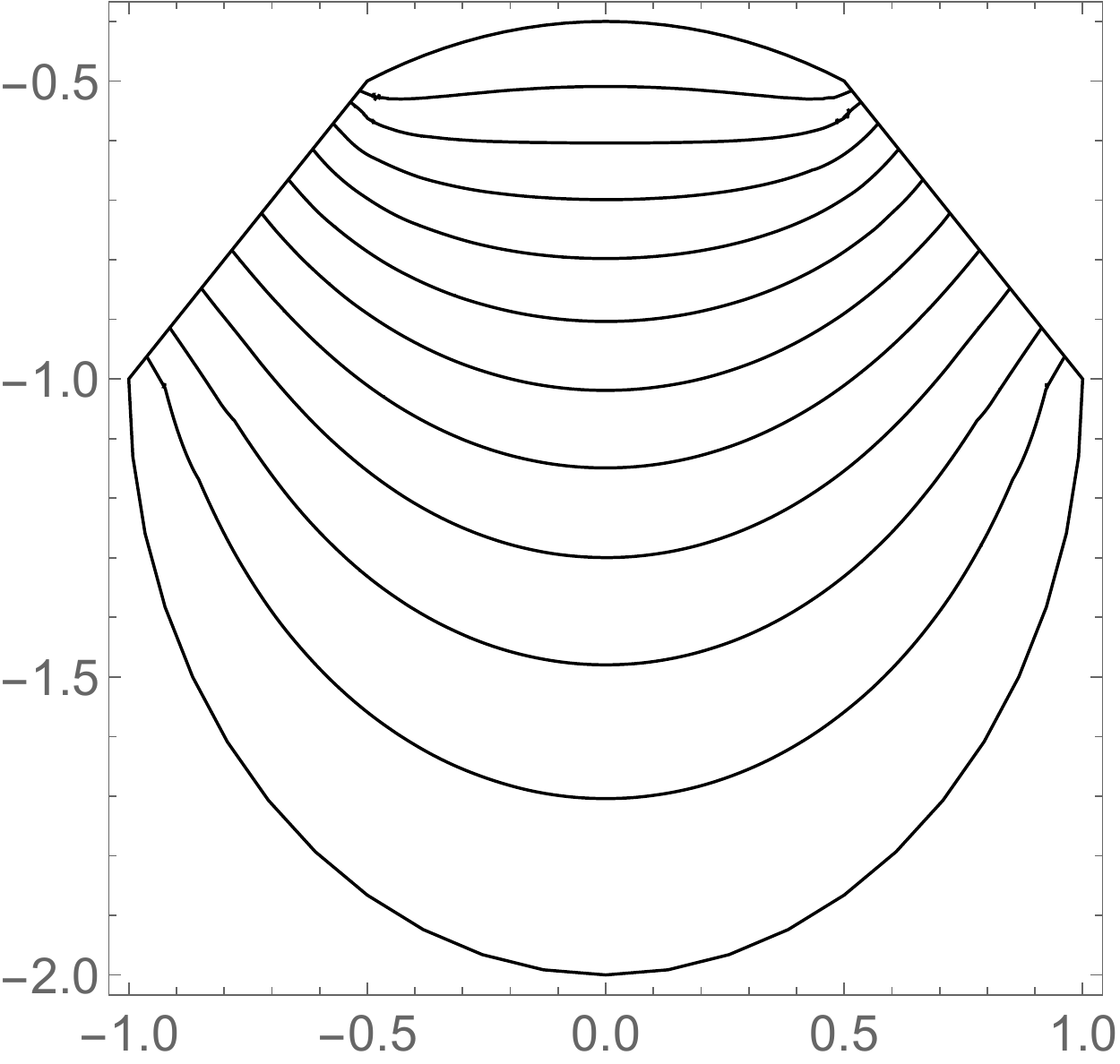}}\quad
    \subfloat[{Potential of the conjugate problem: $\tilde{u}$.}]{\includegraphics[width=0.3\textwidth]{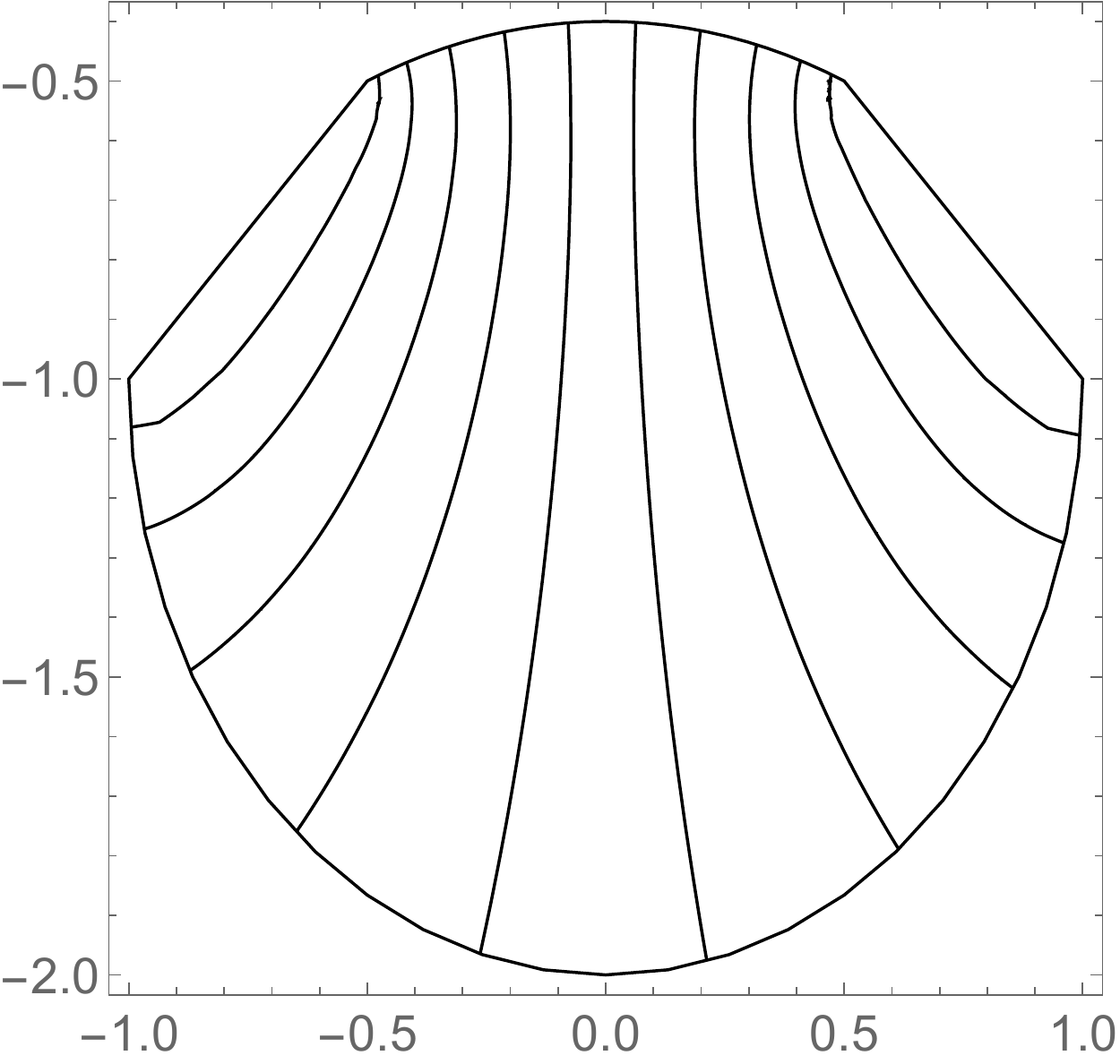}}\quad
    \subfloat[{Map.}]{\label{fig:CGE}\includegraphics[width=0.3\textwidth]{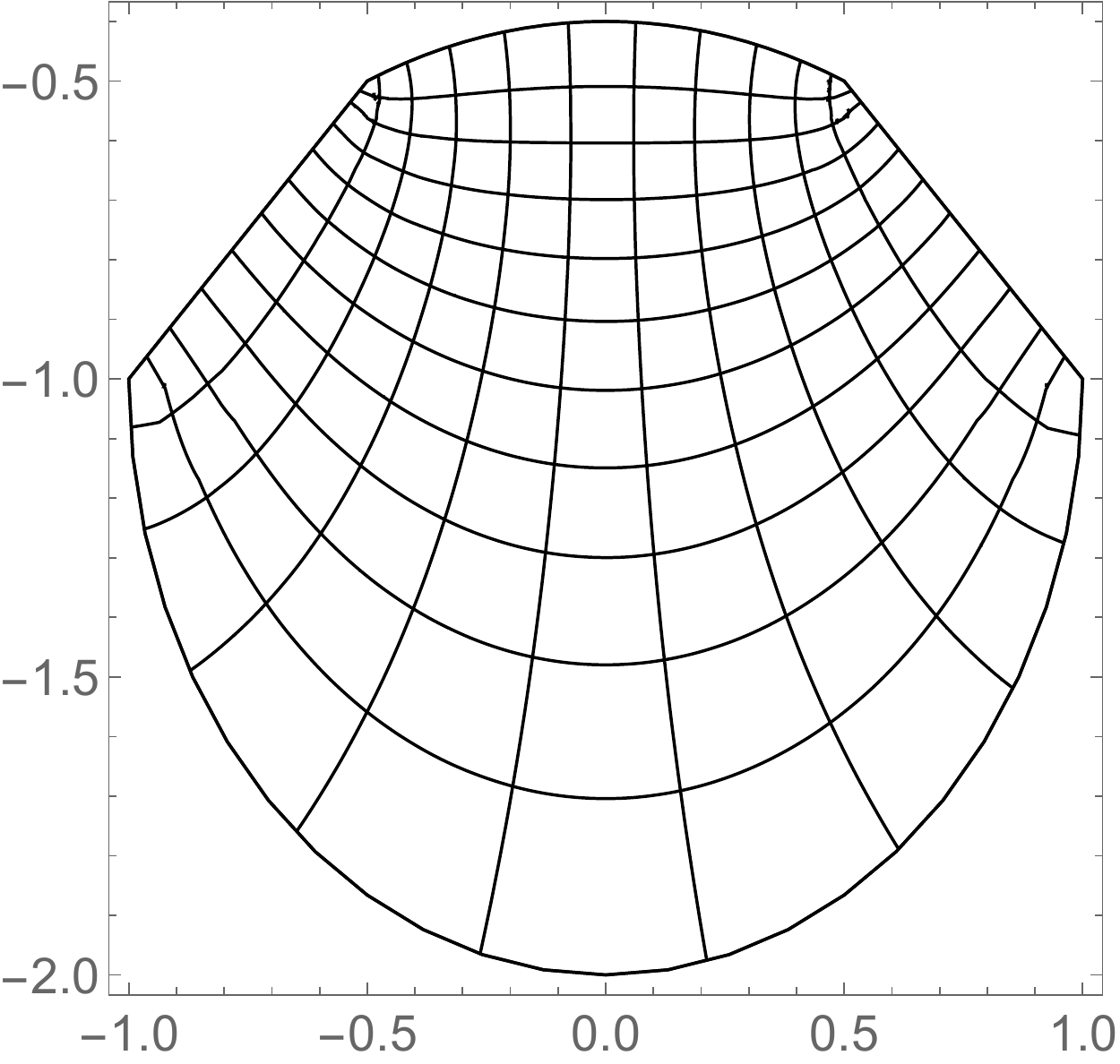}}\quad
    \caption{{Reference capacity via symmetry: $\Omega_{1/4}$. 
    Two error estimates are given in (a) and (b), with error as function of $\sqrt[3]{N}$, where $N$
    is the number of degrees of freedom.
    In both cases the convergence is exponential as predicted by theory.}}\label{fig:CG}
\end{figure}

For the polycircular domain with $m=1$ shown in Figure~\ref{fig:3-cir-m1}, 
even though the exact reference capacity is not known, the domain is such that its inherent symmetries can be exploited.
Consider a quarter domain $\Omega_{1/4}$ formed by connecting the opposing inner and outer arcs with straight edges (see the domain of Figure~\ref{fig:CGC}). Now we define the capacity of the quadrilateral $\Omega_{1/4}$ as in \eqref{eq:u-Dir} except we add zero Neumann conditions on the straight edges.
Using the reciprocal error estimator we can find with high confidence a numerical approximation which is close to machine precision. 
Interestingly, the auxiliary space error estimator agrees with the reciprocal error estimator (see Figures~\ref{fig:CGA} and \ref{fig:CGB}).
For the quarter domain (see Figures~\ref{fig:CGC}-\ref{fig:CGE}) we get 
\begin{equation} \label{eq:refq}
  \capa\,\Omega_{1/4}=1.399386415925543
\end{equation} 
and thus the reference value
\begin{equation} \label{eq:reff}
  \capa\,\Omega = 4\, \capa\,\Omega_{1/4}=5.597545663702171.
\end{equation} 

\begin{figure}
    \centering
    \subfloat[{$m=1$.}]{\label{fig:CGFullA}\includegraphics[width=0.33\textwidth]{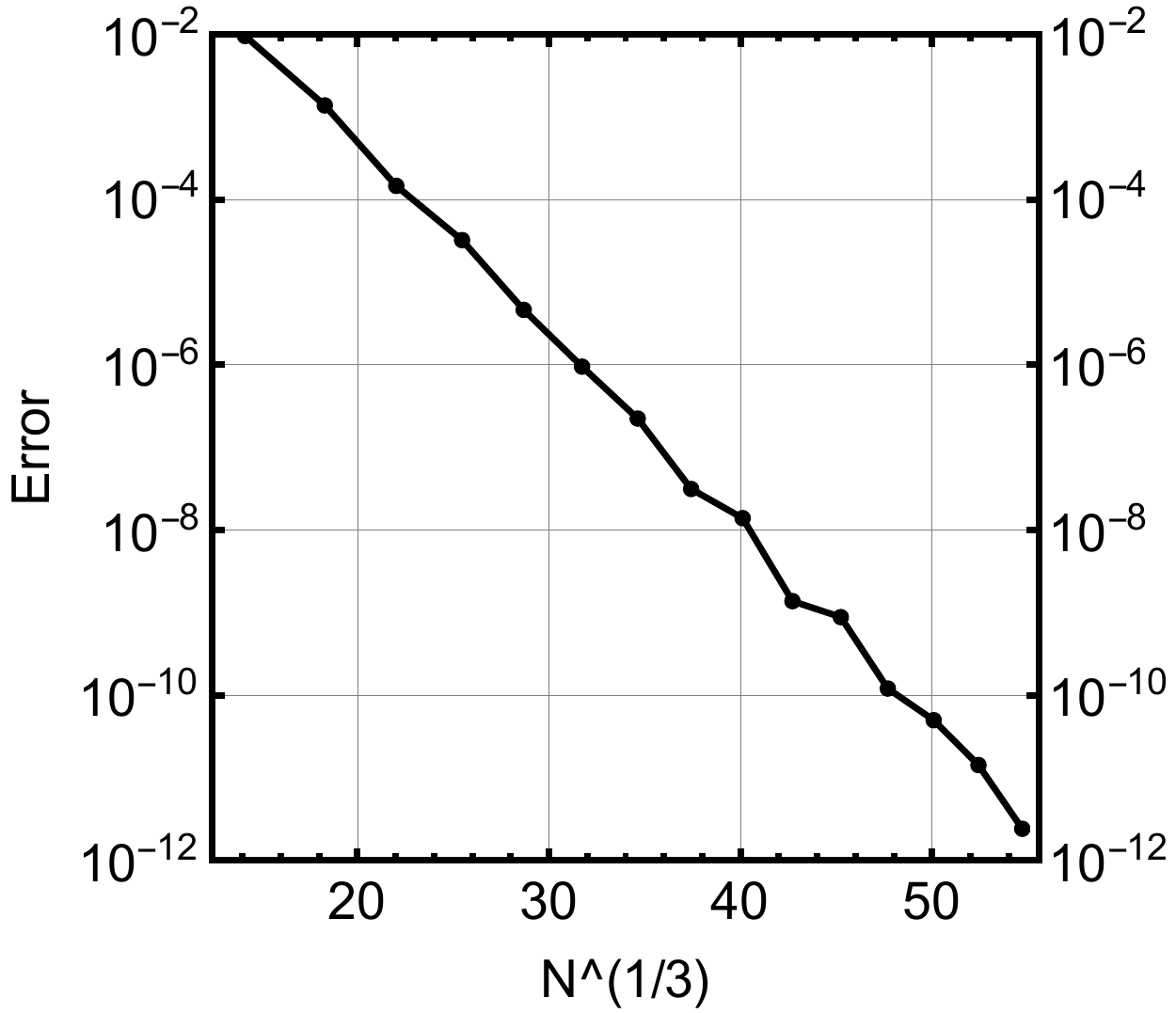}}
    \subfloat[{$m=2$.}]{\label{fig:M2Error}\includegraphics[width=0.33\textwidth]{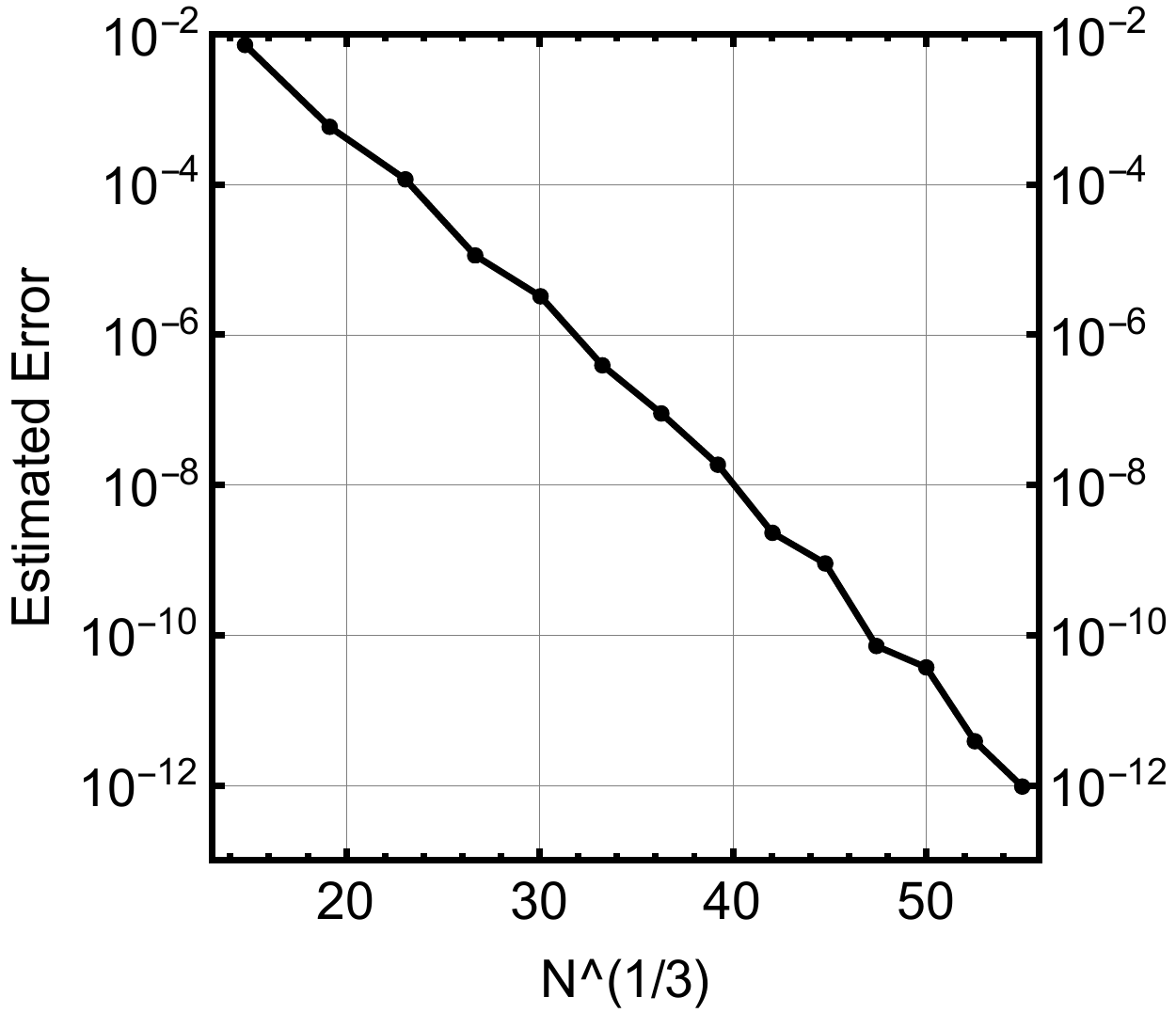}}
    \subfloat[{$m=5$.}]{\label{fig:M5Error}\includegraphics[width=0.33\textwidth]{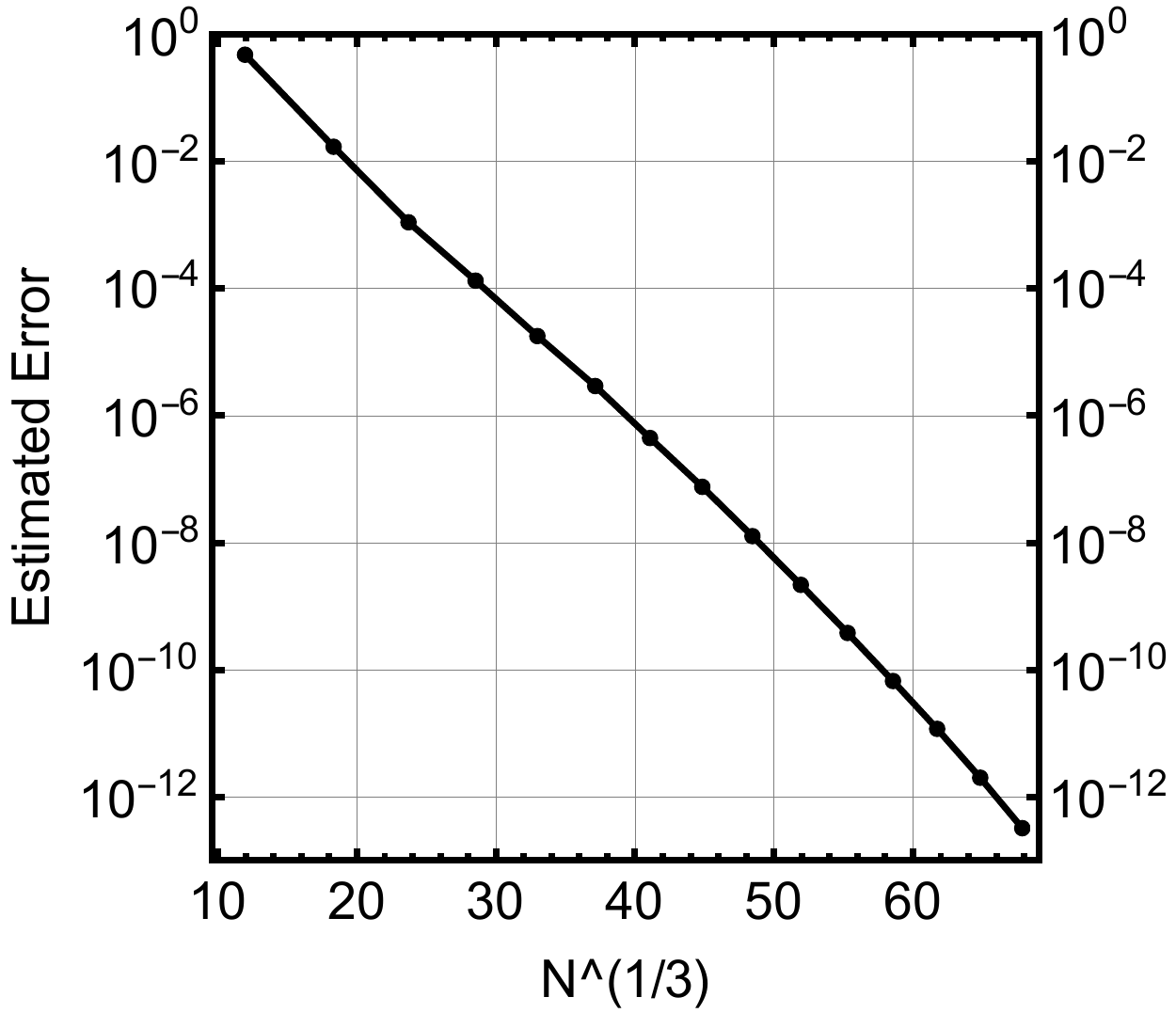}}
    \caption{Error estimates ($hp$-FEM). Logplots with error as function of $\sqrt[3]{N}$. For $m=1$ the reference value
    is computed using symmetries, in the other two cases, the last results in the $hp$-sequence is taken to be the reference.
    In all cases the obtained convergence is exponential with almost identical rate.
    The auxiliary error estimates coincide with the estimated errors and hence are not shown in the graphs.
    The similarity of the graphs gives us confidence in the choice of the reference values for the cases
    $m=2$ and $m=5$.
    }\label{fig:MNError}
\end{figure}

Using the reference value we can estimate the true error of the solution of the full discretization. For the $hp$-FEM, the convergence graph is given in Figure~\ref{fig:CGFullA}. We get
\begin{equation} \label{eq:refC}
  \capa\,\Omega = 5.597545663702324,
\end{equation} 
and observe the error number $= 12$ (13 significant digits). Again, the auxiliary space error estimator
agrees with the effectivity index $\lambda_{m=1} \in [1.14,1.38]$. For the other cases
$\lambda_{m=2} \in [0.91,1.40]$ and $\lambda_{m=5} \in [1.10,1.37]$.

\begin{figure}
    \centering
    \subfloat[{$m=1$.}]{\label{fig:err-BIE1}\includegraphics[width=0.33\textwidth]{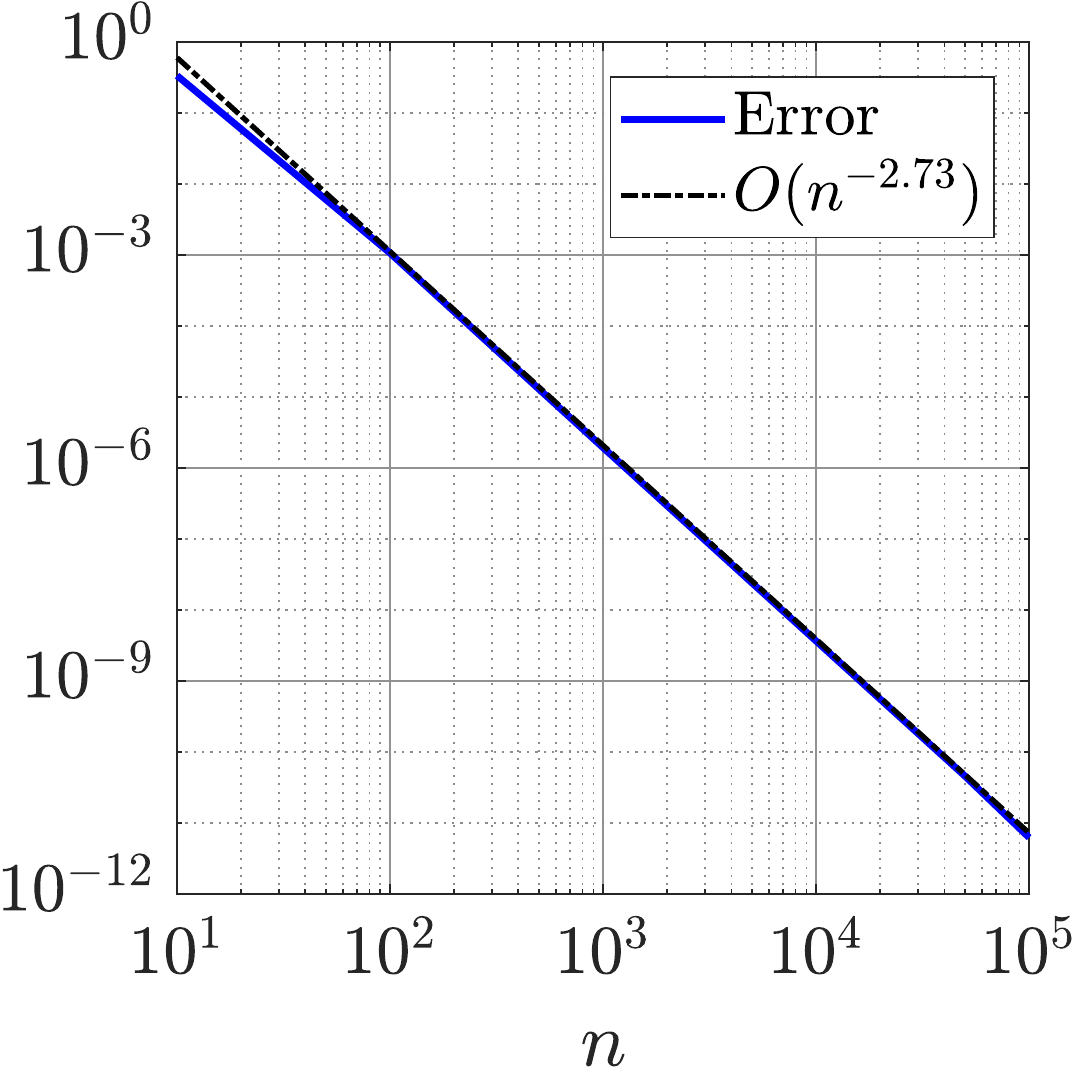}}
    \subfloat[{$m=2$.}]{\label{fig:err-BIE2}\includegraphics[width=0.33\textwidth]{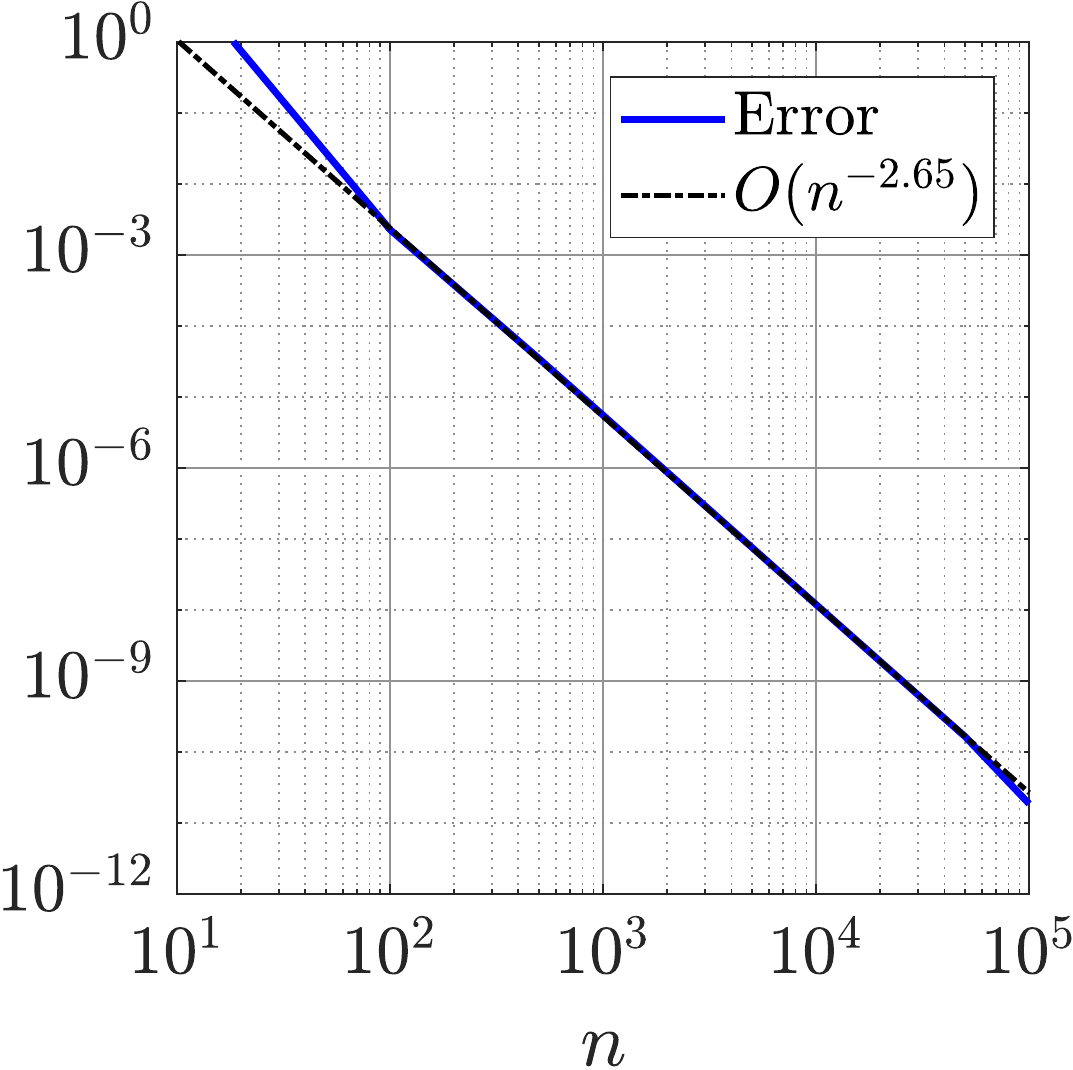}}
    \subfloat[{$m=5$.}]{\label{fig:err-BIE3}\includegraphics[width=0.33\textwidth]{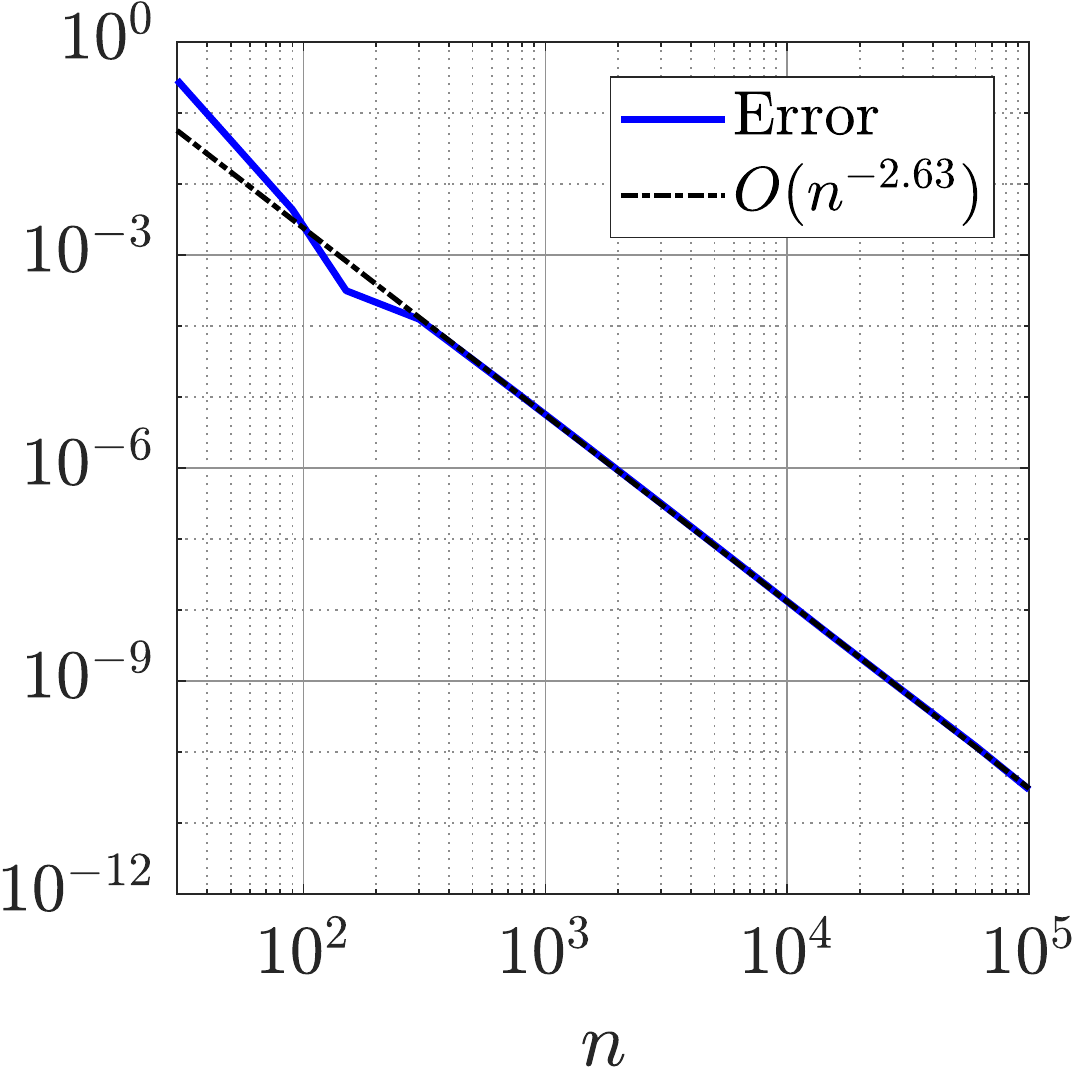}}
    \caption{Error estimates (BIE). The error compared with the reference value in~\eqref{eq:reff} for $m=1$. For $m=2$ and $m=5$, the error is computed by taking the values obtained by $hp$-FEM and presented in Table~\ref{tab:cap-values} as reference values.}\label{fig:err-BIEN}
\end{figure}

For the BIE method, the capacity $\capa\Omega$ is computed using
 the MATLAB function {\tt capm}. We use $n=2^{13}$ and $\alpha=\i$ 
 for $m=1$, $n=2^{13}$ and $\alpha=4\i$ for $m=2$, and $n=15\times2^{9}$ 
 and $\alpha=2\i$ for $m=5$.
For $m=1$, the absolute value of the differences between the approximate 
values and the reference value in~\eqref{eq:reff} vs. $n$ are given in 
Figure~\ref{fig:err-BIE1}. As we can see from Figure~\ref{fig:err-BIE1}, 
the error is $O(n^{-2.73})$. 
For $m=2$ and $m=5$, the absolute error is computed by taking the values 
obtained by $hp$-FEM and presented in Table~\ref{tab:cap-values} as 
reference values. The errors are presented in 
Figures~\ref{fig:err-BIE2}--\ref{fig:err-BIE3} which are
 $O(n^{-2.65})$ and $O(n^{-2.63})$, respectively.
The presented results illustrate that the BIE method achieves a comparable level of accuracy with the $hp$-FEM.

\begin{figure}
    \centering
    \includegraphics[width=0.65\textwidth]{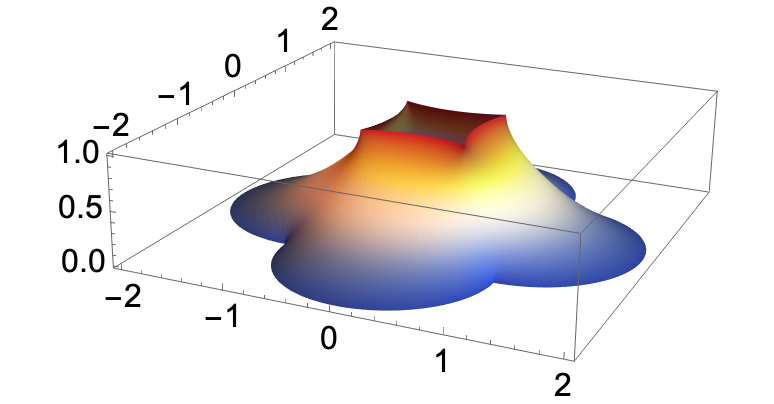}
    \caption{{Potential: $m=1$.}}\label{fig:CGFull}
\end{figure}

%FILE: sec511.tex
%%%%%%%%%%%%%%%%%%%%%%%%%%%%%%%%%%%%%%%%%
\section{Numerical experiments on domains with challenging boundaries}

In this section we study three different aspects of polycircular domains.
First, in the absence of exact capacity values, we study M\"obius invariance
of a given configuration. Here again, we can employ symmetries and obtain a numerical
reference value of high accuracy. Second, we consider a domain with a jagged boundary
with small angles, an interesting test case for BIE. 
Finally, a parameter-dependent configuration
with lens-shaped boundary is tested against upper and lower bounds given in the literature.
Also in this case symmetries can be exploited.

\begin{figure}[H]
  \centerline{
	\includegraphics[width=0.3\linewidth]{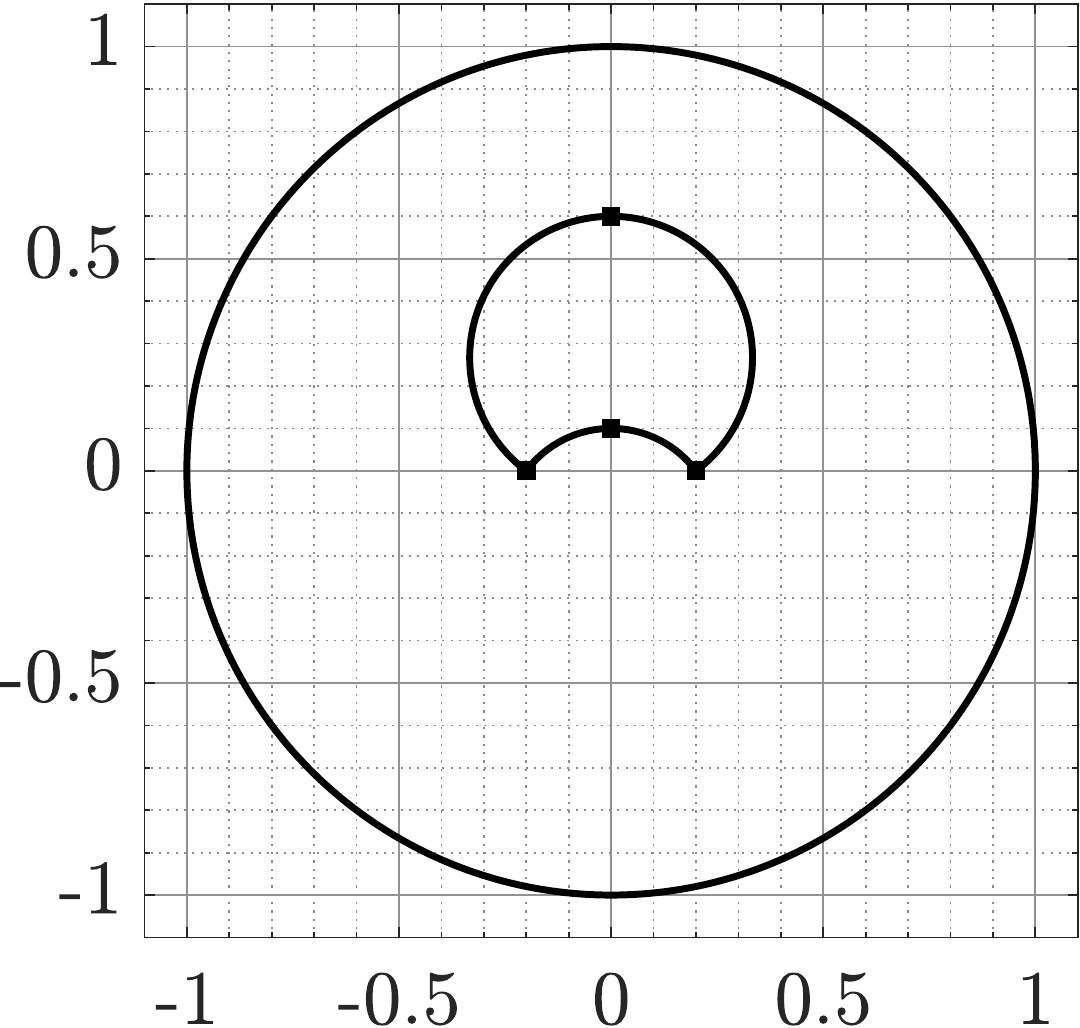}
	\hfill
	\includegraphics[width=0.3\linewidth]{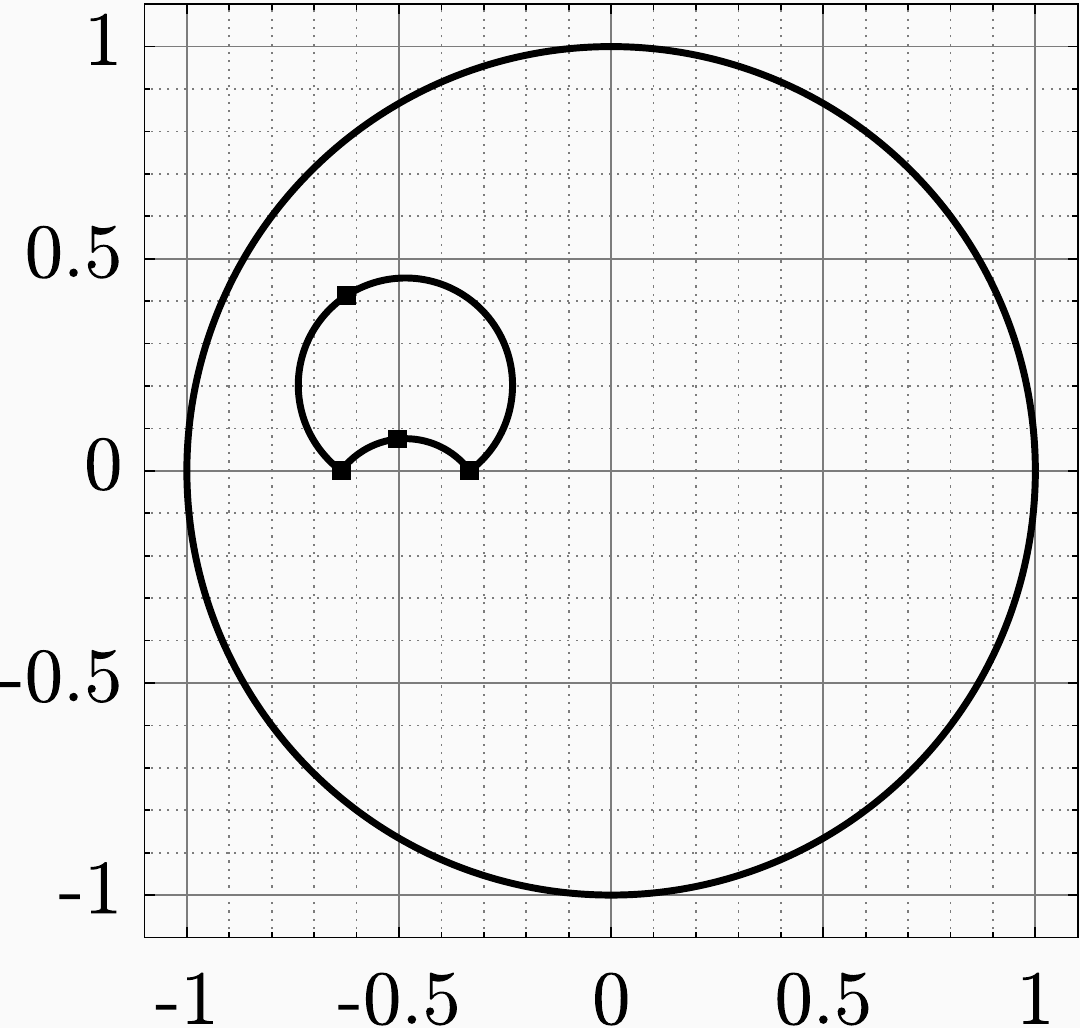}
	\hfill
	\includegraphics[width=0.3\linewidth]{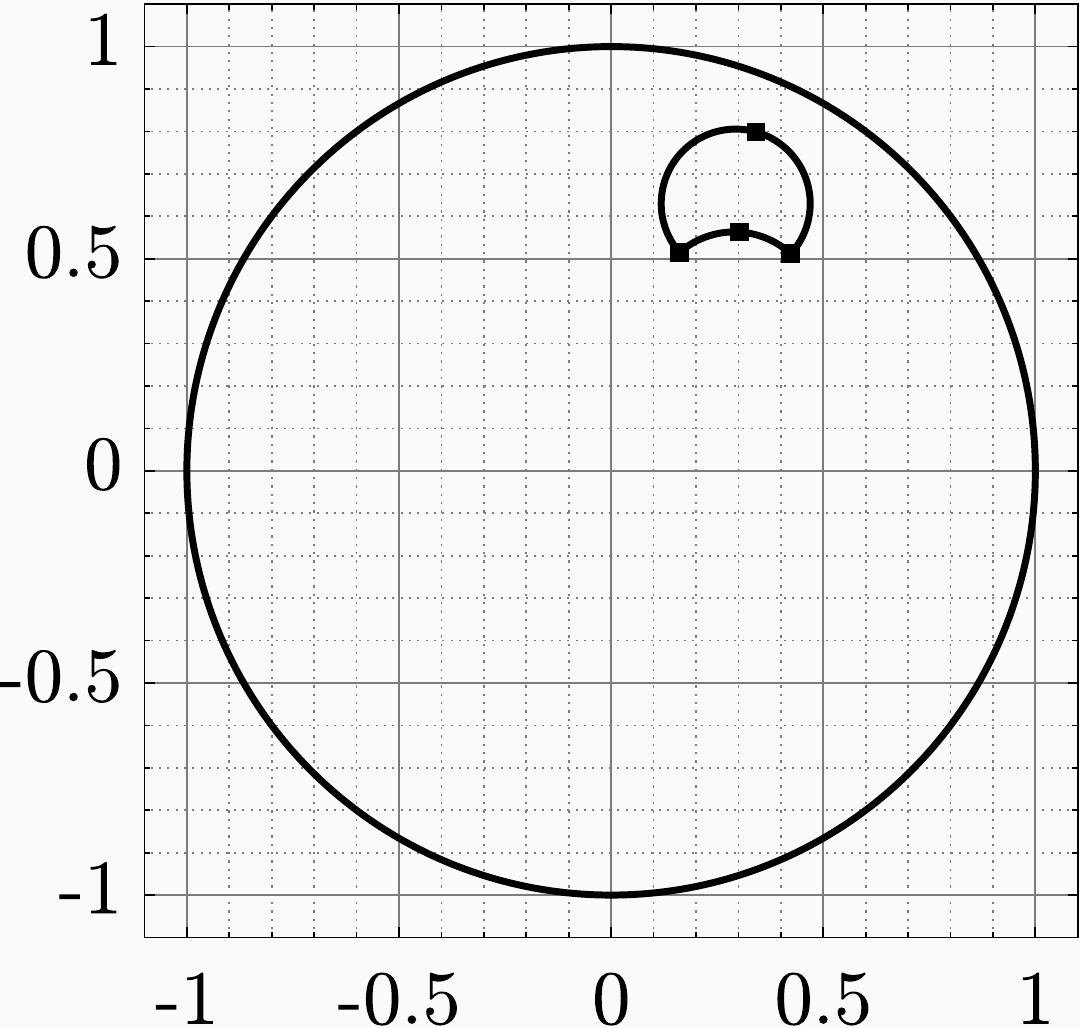}
	}
  \caption{The condenser $(\B^2,E)$ (left) and the condenser $(\B^2,T_a(E))$ for $a=0.5$ (center), and $a=-0.3-0.5\i$ (right).}
  \label{fig:Mob}
\end{figure}

%%%%%%%%%%%%%%%%%%%%%%%%%%%%%%%%%%%%%%%%%
\subsectionB{\textit{Checking the  M\"obius invariance}\\}
%\subsection{Checking the  M\"obius invariance}
%In this section, we study several condensers $(\B^2,E)$ where $\B^2$ is the unit disk and $E$ is a connected circular-arc polygon. 
Unfortunately finding examples of  polycircular   condensers with 
known capacity looks too difficult. Therefore, one way to check the accuracy 
of the presented methods is to use the fact that capacity is conformal invariant. 
Consider the condenser $(\B^2,E)$ where $E$ is the closure of the domain bordered by the circular arc polygon consists of two circular arcs, the first one passes through the three points $-2$, $0.6\i$, and $0.2$, and the second one passes through the points $0.2$, $0.1\i$, and 
$-0.2$ (see Figure~\ref{fig:Mob}). Let 
\[
T_a(z) = \frac{z-a}{1-\overline{a} z}
\]
be a M\"obius transformation where $a\in \B^2$. Then, by the conformal invariance of the capacity, 
we should have
\[
\capa(\B^2,E) = \capa(\B^2,T_a(E)).
\]

The two numerical methods are used to compute $\capa(\B^2,T_a(E))$ for several values of $a$ and 
\[
\left|\capa(\B^2,E) - \capa(\B^2,T_a(E))\right|
\]
is considered to be the error in the computed values. The obtained results are presented in Table~\ref{tab:Mob} for several values of $a$ where $T_0(E)=E$. For the BIE method, we use $n=2^{12}$, $\alpha=T_a(-0.3\i)$, and $\alpha_1=T_a(0.3\i)$. Both methods, the BIE method and the $hp$-FEM, give results with almost the same level of high accuracy.

\begin{table}[H]
\caption{The approximate values $\capa(\B^2,T_a(E))$ and the error $\left|\capa(\B^2,E) - \capa(\B^2,T_a(E))\right|$ where $\capa(\B^2,E)=\capa(\B^2,T_0(E))$. The two methods agree with 10 significant digits.
The reference value for $a=0$ is $6.044918141954128$.}\label{tab:Mob}
\centering
\begin{tabular}{|l|c|c||c|c|}
    \hline
   $a$ & \multicolumn{2}{c||}{BIE Method} & \multicolumn{2}{c|}{$hp$-FEM} \\
    \cline{2-5}
	   & $\capa(\B^2,T_a(E))$ & Error & $\capa(\B^2,T_a(E))$ & Error \\ \hline
 $0$          & 6.04491814238562 & --- & 6.044918141954396 & --- \\
 $0.1$        & 6.04491814238563 & $9.77\times10^{-15}$ & 6.044918141954359 & $3.73\times10^{-14}$\\
 $0.5$        & 6.04491814238567 & $4.88\times10^{-14}$ & 6.044918141954392 & $3.55\times10^{-15}$\\
 $0.1+0.3\i$  & 6.04491814238559 & $2.93\times10^{-14}$ & 6.044918141954345 & $5.06\times10^{-14}$\\
 $-0.2+0.5\i$ & 6.04491814238574 & $1.18\times10^{-13}$ & 6.044918141954328 & $6.75\times10^{-14}$\\
 $-0.3-0.5\i$ & 6.04491814238566 & $3.55\times10^{-14}$ & 6.044918141954332 & $6.39\times10^{-14}$\\
    \hline
\end{tabular}
\end{table}

%\end{mysubsection}

\begin{figure}
\centering
\begin{tabular}{p{7cm} p{7cm}}
    \vspace{0pt} 
    \subfloat[The domain.]{\label{fig:BSa}\includegraphics[width=7cm]{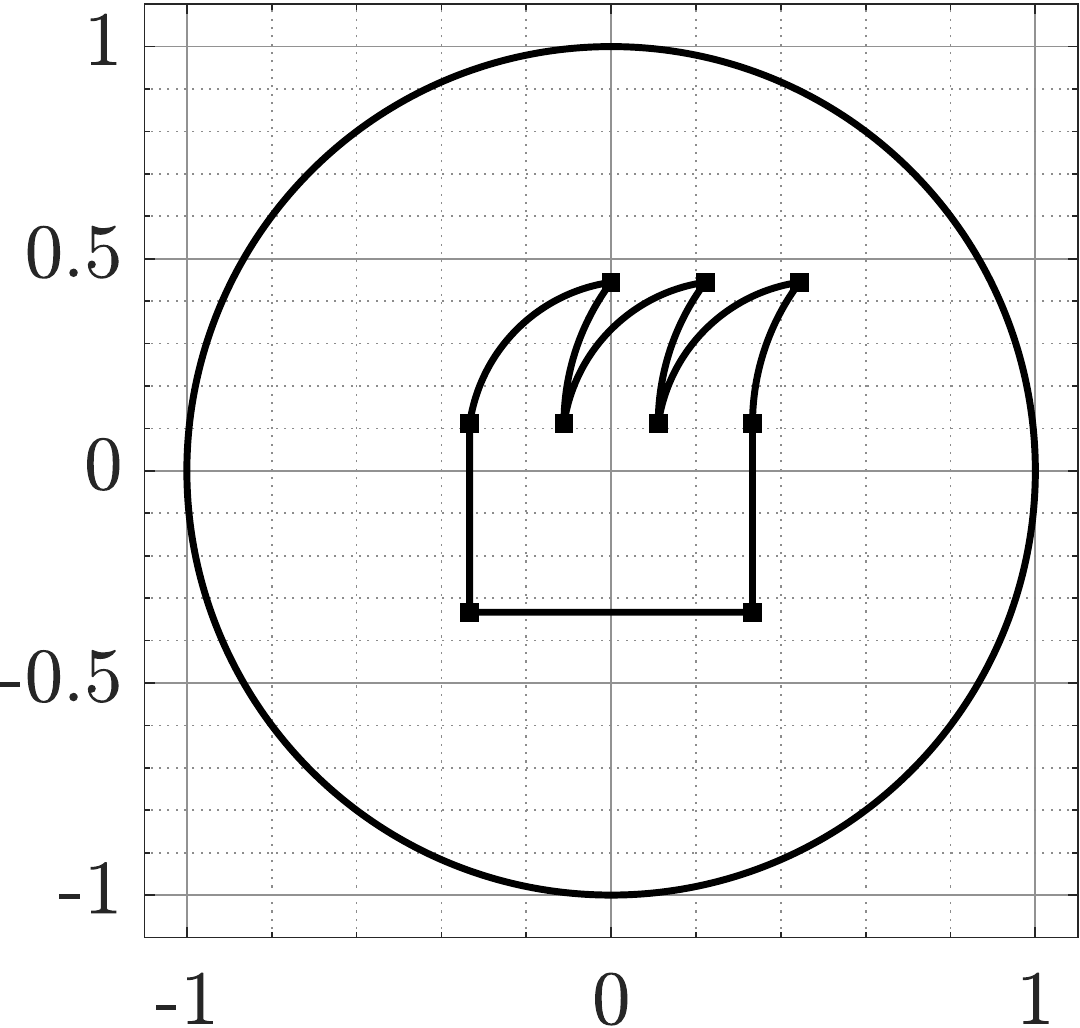}}
    & 
    \vspace{0pt}  \subfloat[BIE convergence.]{\label{fig:BSb}
\begin{tabular}{l|l}
    \hline
    $n$ & $\capa(\B^2,E)$ \\ \hline
    $9\times2^8$     & $7.26568934865634$ \\
    $9\times2^9$     & $7.26568351878449$ \\
    $9\times2^{10}$  & $7.26568260732833$ \\
    $9\times2^{11}$	 & $7.26568246621964$ \\ \hline
\end{tabular}
}
\end{tabular}

  \caption{Bart Simpson condenser. The domain and capacities computed with BIE.}
  \label{fig:BS}
\end{figure}

%%%%%%%%%%%%%%%%%%%%%%%%%%%%%%%%%%%%%%%%%
\subsectionB{\textit{Bart Simpson condenser}\\}
%\subsection{Bart Simpson condenser}
We assume that $E$ has a Bart Simpson shape as shown in
Figure~\ref{fig:BS}. The vertices of the $E$ are $3/9+\i/9$,
$3/9-3\i/9$, $-3/9-3\i/9$, $-3/9+\i/9$, $4\i/9$, $-1/9+\i/9$,
$2/9+4\i/9$, $1/9+\i/9$, and $4/9+4\i/9$. The arcs in
Figure~\ref{fig:BS} are determined either by their end-points and a point on
each arc, which are $-2/9+3\i/9$, $3\i/9$, $2/9+3\i/9$, respectively, or 
by their end-points and the center of the circle containing
the arc, which are $4/9+\i/9$, $6/9+\i/9$, and $8/9+\i/9$, respectively. 
Note that the inner boundary has no cusps as the outer angles at the vertices $-1/9+\i/9$ and $1/9+\i/9$ are approximately $8^{\rm o}$. 

Here the capacity obtained with $hp$-FEM with $p=18$ is $\capa(\B^2,E) = 7.265682491066263$ and the corresponding error estimate is $\sim 10^{-8}$.
The BIE method used with $n=9\times2^{11}$, $\alpha=-0.5\i$, and $\alpha_1=0$ gives the approximate value 
$\capa(\B^2,E)=7.26568246621964$ which agrees with eight significant digits with the approximate value obtained by the FEM method. 
The difference between the two results is within the auxiliary space error estimate.
Further, the approximate values for $\capa(\B^2,E)$ obtained for several values of $n$ are given in the table in Figure~\ref{fig:BSb}.
%\begin{center}
%\begin{tabular}{l|l|l|l|l}
%    \hline
%$n$             & $n=9\times2^8$     & $n=9\times2^9$     & $n=9\times2^{10}$  & $n=9\times2^{11}$ 
%\\ \hline
%$\capa(\B^2,E)$ & $7.26568934865634$ & $7.26568351878449$ & $7.26568260732833$ & $7.26568246621964$\\
%\hline
%\end{tabular}
%\end{center}
These presented results demonstrate that the BIE method gives accurate results even though the angles of two corner points are small (approximately equal $8^{\rm o}$).

%    \begin{equation}
% {\rm H. Hakula:\,\, cap}  = 
% 7.265682491066263
%     ; \quad 2021-08-05
%    \end{equation}   

%\end{mysubsection}

\begin{figure}
\centering
  \subfloat[The domain $\DD\backslash E$ for the lens $E$.]{
	\includegraphics[width=0.45\linewidth]{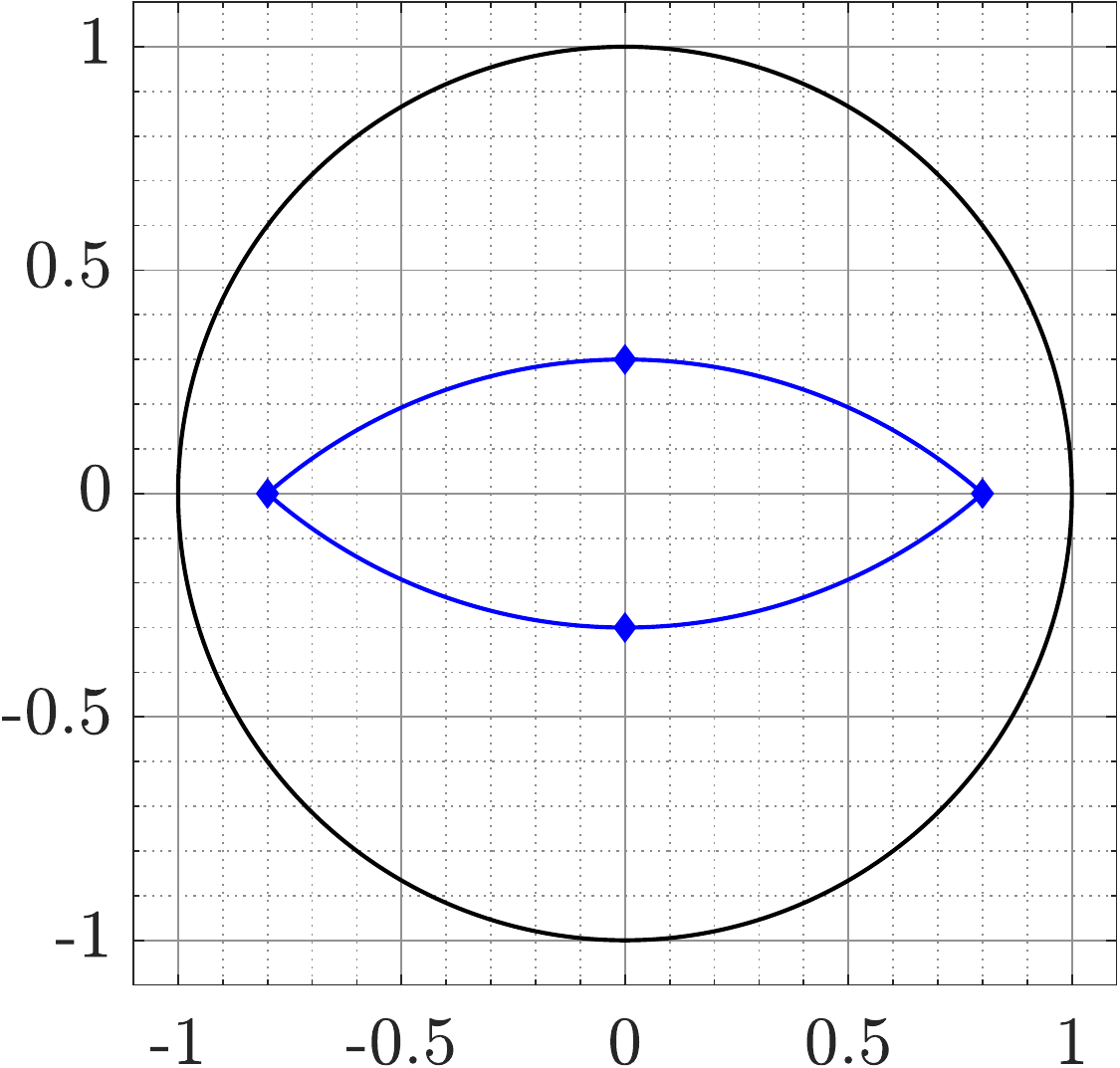}
	}\quad
  \subfloat[Convergence of BIE.]{
	\includegraphics[width=0.45\linewidth]{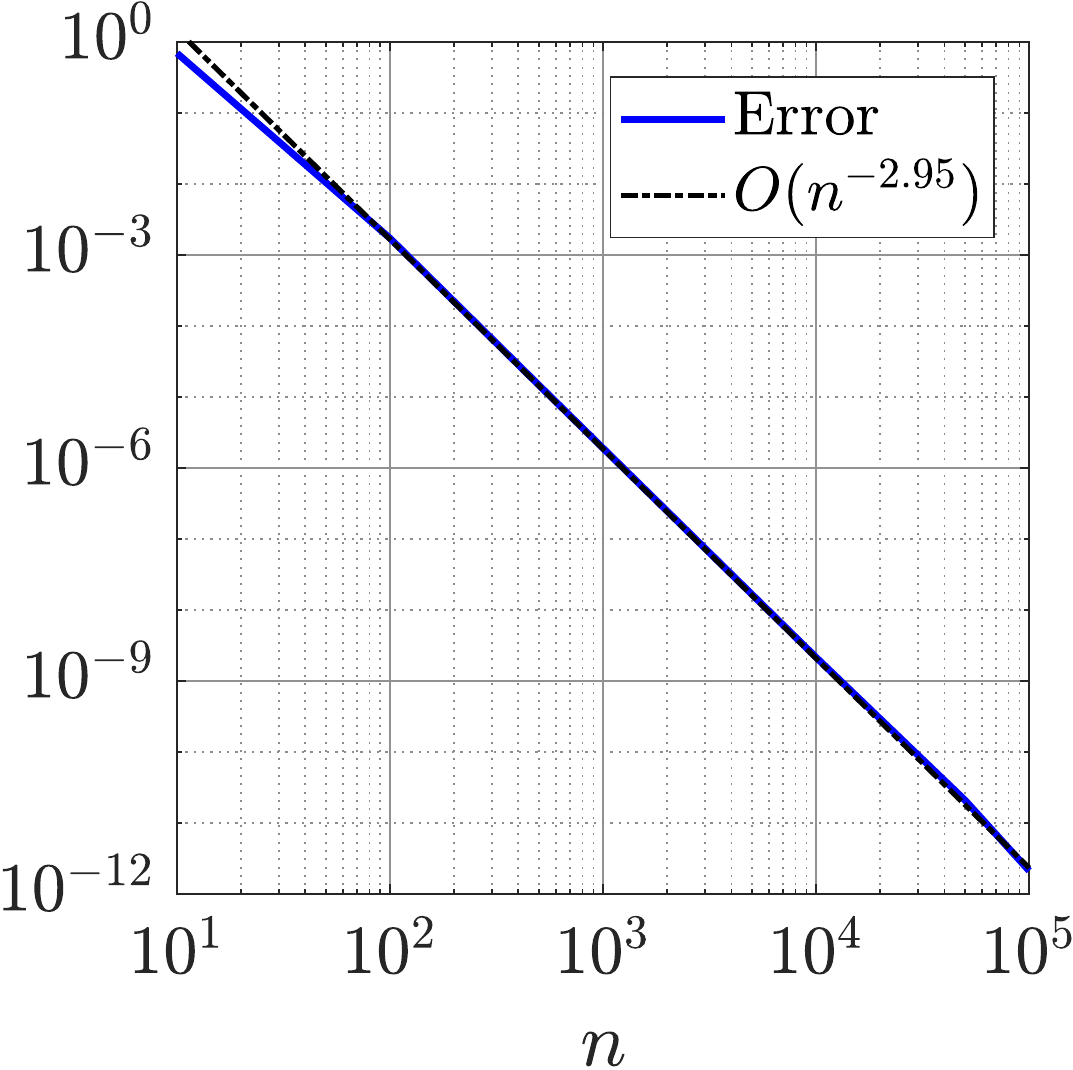}
	}
	\caption{Lens domain with $r=4/5$, $s=3/10$. The reference value is $10.15585205509004$.}
  \label{fig:lens-1}
\end{figure}

\begin{figure}
\centering
  \subfloat[Capacity as function of $s$.]{
	\includegraphics[width=0.45\linewidth]{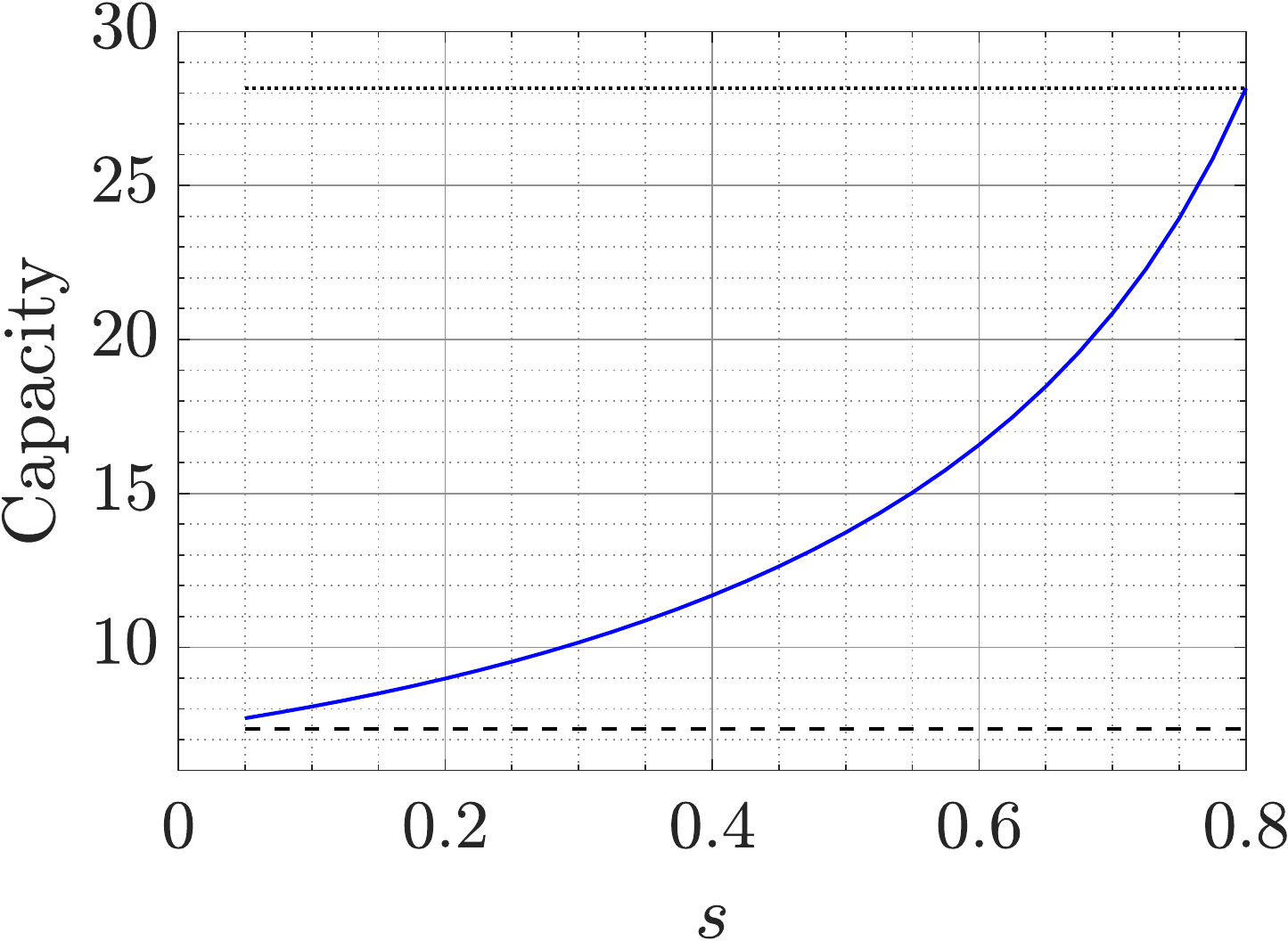}
	}\quad
  \subfloat[Capacity as function of hyperbolic perimeter.]{
	\includegraphics[width=0.45\linewidth]{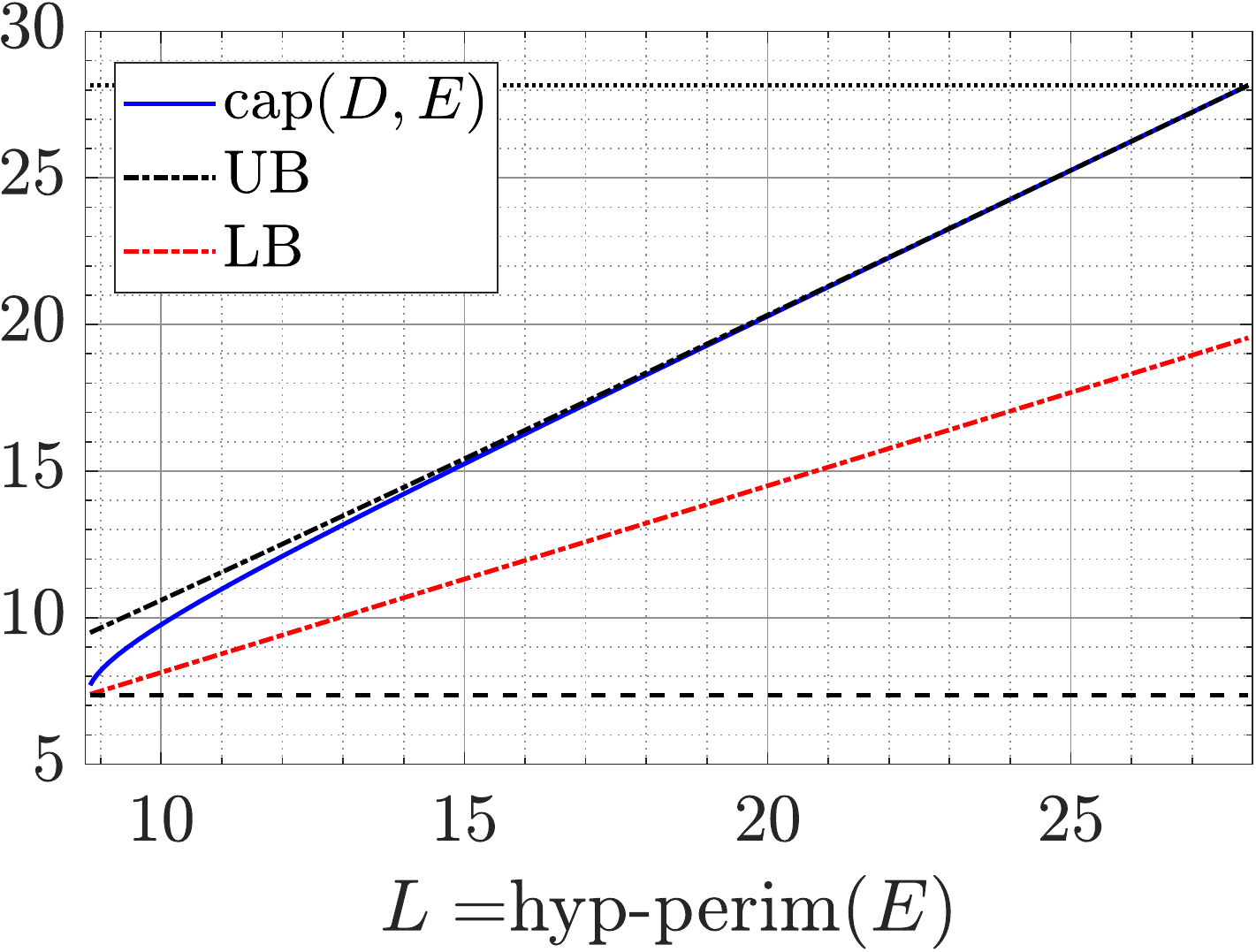}
	}
	\caption{The capacity $\capa(\DD,E)$, the upper bound~\eqref{eq:UB}, and the lower bound~\eqref{eq:LB} where $E$ is the lens domain with $r=0.8$ and $0<s\le r$ vs. the hyperbolic perimeter of $E$. The horizontal dot line is $\capa(\DD,\overline{B(0,r)})$ and the horizontal dashed line is $\capa(\DD,[-r,r])$.}
  \label{fig:lens-2hp}
\end{figure}

%\begin{mysubsection}{\bf Lens domains.}
\subsectionB{\textit{Lens domains}\\}
%\subsection{Lens domains}
Consider the condenser $(\DD,E)$ where $E$ is the closure of the lens domain bordered by the 
two circular arcs, the first one passes through the points $r,-\i s,-r$, and 
the second one passes through the points $-r, \i s,r$ 
(see Figure~\ref{fig:lens-1}a for $r=0.8$ and $s=0.3$).
 The BIE method is used to compute the capacity of the condenser 
 $(\DD,E)$ for $r=0.8$ and for $0.05\le s\le r$. The obtained results are shown 
 in Figure~\ref{fig:lens-2hp}. 
 The convergence graph for a single instance $r=0.8$ and $s=0.3$ is shown in Figure~\ref{fig:lens-1}b.
 The values of the capacity $\capa(\DD,E)$ vs. the hyperbolic perimeter of $E$, $L=\mbox{hyp-perim}(E)$, are presented in Figure~\ref{fig:lens-2hp}. It is known that~\cite[Eqs~(5.13)--(5.15)]{nrv2} 
\begin{equation}\label{eq:UB}
\frac{2\pi}{\log\left(\sqrt{1+(2\pi/L)^2}+2\pi/L\right)}
\end{equation}
is an upper bound for $\capa(\DD,E)$ and, since $E$ is convex, 
\begin{equation}\label{eq:LB}
\frac{2\pi}{\mu\left(\th(L/4)\right)}
\end{equation}
is a lower bound for $\capa(\DD,E)$, where~\cite[p.~122]{HKV}
\begin{equation} \label{capGro}
\mu(r)=\frac{\pi}{2}\frac{\K'(r)}{\K(r)}, \quad \K(r)=\int^1_0 \frac{dx}{\sqrt{(1-x^2)(1-r^2x^2)}},
\quad \K'(r)=\K(r'), \quad r'=\sqrt{1-r^2}\,.
\end{equation}
Here, $\K(r)$ and $\K'(r)$ are the complete elliptic integrals of the first kind.
These bounds are shown in Figure~\ref{fig:lens-2hp} as well.

Note that for $s=r=0.8$, the set $E$ is the closed disk $\overline{B(0,r)}=\{z:|z|\le r\}$ and the capacity $\capa(\DD,E)$ for the latter case is 
 $2\pi/\log(1/r)$ and hence, for $r=0.8$, equal to $28.157593038985901$ 
 which is the horizontal dotted-line in Figure~\ref{fig:lens-2hp}. 
 In the case $s=0$, the set $E=[-r,r]$, and hence, by the M\"obius  transform 
\[
z\mapsto \frac{z+r}{1+rz},
\]
the domain $\DD\backslash E$ is mapped to the domain $\DD\backslash \hat E$
 where $\hat E=[0,2r/(1+r^2)]$. Hence, 
\[
\capa(\DD,E) = \capa(\DD,\hat E) =\frac{2\pi}{\mu\left(\frac{2r}{1+r^2}\right)}
\]
and we see that for $r=0.8$, $\capa(\DD,E)=7.360222723821019,$ which value is drawn as the 
horizontal dashed-line in Figure~\ref{fig:lens-2hp}.

\begin{figure}[H]
\begin{tabular}{p{7cm} p{8cm}}
    \vspace{0pt} 

  \subfloat[The domain.]{
	\includegraphics[width=6cm]{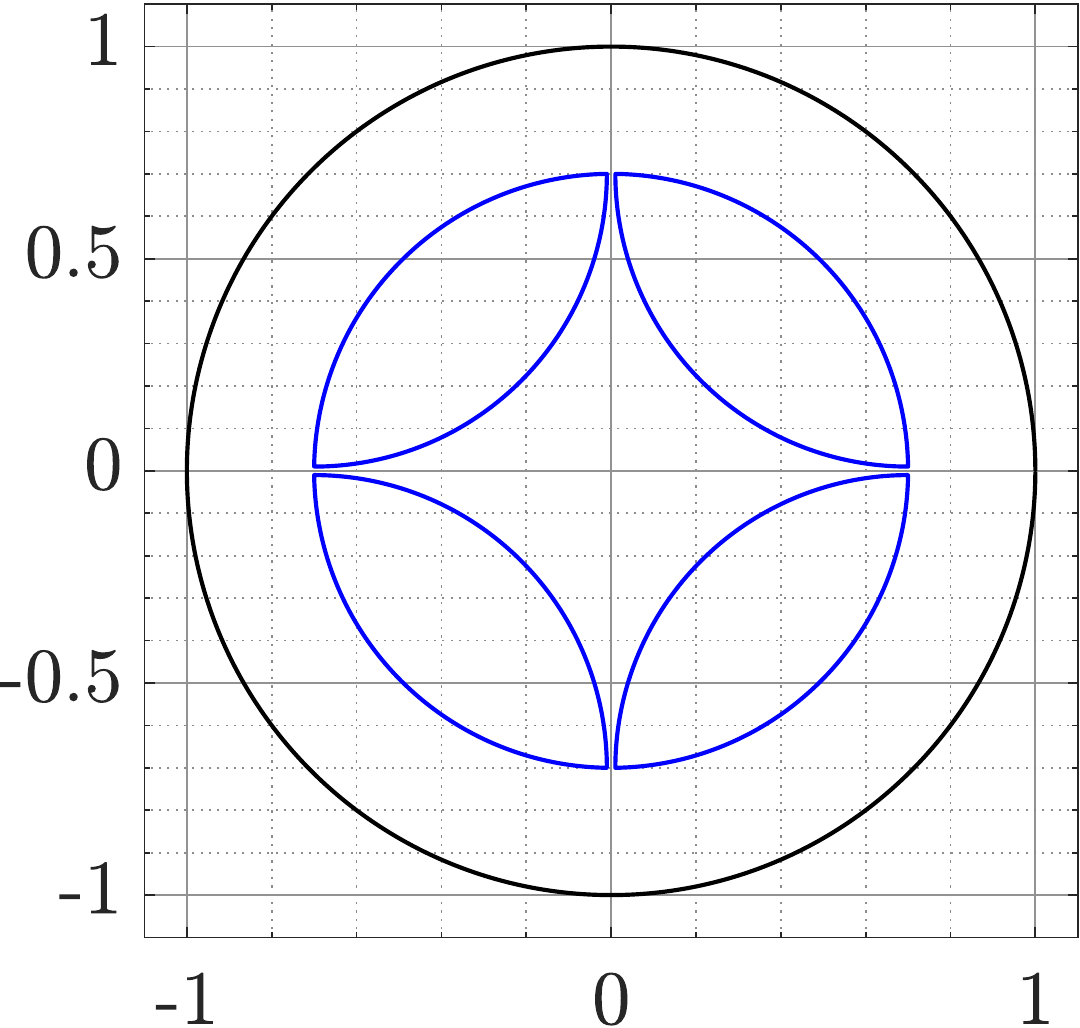}
	}
	&
	\vspace{0pt}  
	\subfloat[Single lens configuration.]{\label{tab:lens-2}
\begin{tabular}{l|l|l}
    \hline
 Lens 1    & Arc 1                 & Arc 2           \\ \hline
 Point 1   & $0.01+0.7\i$          & $0.7+0.01\i$   \\
 Point 2   & $0.7(1+\i)/\sqrt{2}$  & $0.3(1+\i)/\sqrt{2}$\\     
 Point 3   & $0.7+0.01\i$          &$0.01+0.7\i$  \\ \hline
\end{tabular}
	}
\end{tabular}
	\caption{The domain $\DD\backslash E$ where $E$ consists of four lenses.
	The four lenses are defined as rotated images of the configuration given in (b).
	}
  \label{fig:lens-2}
\end{figure}

For the second example, we consider the condenser $(\DD,E)$ where $E$ consists of 
four lens domains as in Figure~\ref{fig:lens-2}. Each of these lenses is defined
 by two circular arcs, each of which passes through three given points. 
 The values of these points are given in Table~\ref{tab:lens-2}. The computed 
 value of the capacity is 
 
\begin{equation} \label{caplens1}
 \capa(\DD,E)=17.613666396960355.
\end{equation}

 Note that this value of 
 the capacity is close to the value of the capacity  of the annulus $0.7<|z|<1$ 
 which is $17.615998583457760$. 

%\begin{table}
%\caption{The points for the lens domain in Figure~\ref{fig:lens-2}.}\label{tab:lens-2}
%\centering
%\begin{tabular}{l|l|l|l|l}
%    \hline
%         &         & Point 1       & Point 2               & Point 3 \\ \hline
%Lens 1   & Arc 1   & $0.01+0.7\i$  & $0.7(1+\i)/\sqrt{2}$  & $0.7+0.01\i$\\
%         & Arc 2   & $0.7+0.01\i$  & $0.3(1+\i)/\sqrt{2}$  & $0.01+0.7\i$\\     \hline
%Lens 2   & Arc 1   & $-0.7+0.01\i$ & $0.7(-1+\i)/\sqrt{2}$ & $-0.01+0.7\i$ \\
%         & Arc 2   & $-0.01+0.7\i$ & $0.3(-1+\i)/\sqrt{2}$ & $-0.7+0.01\i$ \\     \hline
%Lens 3   & Arc 1   & $-0.01-0.7\i$ & $0.7(-1-\i)/\sqrt{2}$ & $-0.7-0.01\i$ \\
%         & Arc 2   & $-0.7-0.01\i$ & $0.3(-1-\i)/\sqrt{2}$ & $-0.01-0.7\i$ \\     \hline
%Lens 4   & Arc 1   & $0.7-0.01\i$  & $0.7(1-\i)/\sqrt{2}$  & $0.01-0.7\i$\\
%         & Arc 2   & $0.01-0.7\i$  & $0.3(1-\i)/\sqrt{2}$  & $0.7-0.01\i$\\     \hline
%\end{tabular}
%\end{table}

%\end{mysubsection}

%%%%%%%%%%%%%%%%%%%%%%%%%%%%%%%%%%%%
\section{Conclusion}

{Polycircular  domains are domains whose boundaries are 
 unions of arcs of circles and rectilinear segments.}
%General polycircular domains can be multiply connected and contain cusps leading to strong singularities.
Computation of conformal capacity of such domains in closed form is possible
only in simple geometric configurations, and numerical approximations are 
necessary. As always with approximations, the central question is the 
accuracy of the obtained results. Within the PDE context the existing 
methods such as $hp$-FEM are equipped with error estimators,
that can be shown to be asymptotically correct. 
Here we employ another approach for verification -- we apply 
two different numerical methods, $hp$-FEM and the boundary integral method, with different discretization characteristics.
Theoretical error bounds are not yet available in the case of the method based on the BIE with the generalized Neumann kernel.
On the other hand, the BIE method
has turned out to be very flexible and easily adjustable to great many geometries
with high precision in the cases reported in \cite{nv1,nv}. 
The numerical experiments presented in this paper gave further evidence about the accuracy of this BIE method and illustrate that the BIE method gives results with a comparable level of accuracy with the $hp$-FEM.
Furthermore, the BIE method is not only easier to implement but also consistently faster both in terms of problem setup and solution time within the class of problems considered here. For an example implementation, we refer to \url{https://github.com/mmsnasser/polycircular}.

Over a carefully designed set of numerical experiments we first establish the convergence properties of the methods, then verify the agreement against known reference results, and finally compare the results.
The computed capacities are extremely accurate and can serve as reference values. In some cases even practical machine precision is obtained. 
More importantly, in those cases where the results differ, though never over the single precision, the
discrepancy between the methods is within the asymptotic error estimates. This gives us very high confidence that the numbers reported here are indeed as accurate as we claim.

In conclusion, we have used two methods to compute capacity of condensers which requires solving a Dirichlet boundary value problem in multiply connected domains for the harmonic potential function whose values are equal to zero on the outer boundary component and equal to one on all inner boundary components. 
The presented methods can be used to compute other quantities in potential theory which require computing harmonic functions with other type of boundary conditions, see~\cite{nv1}, \cite{hare}. 
For example, recently, an analytic formula for computing a quantity known as the ``h-function'' for multiply connected slit domains is presented in~\cite{Green}. The numerical computation of the ``h-function'' for multiply connected slit domains can be formulated as a boundary value problem of the type considered in~\cite{nv1} and hence the $h$-function can be computed for multiply connected domains with high connectivity by the BIE method.

\section*{Acknowledgments}
The authors would like to thank two anonymous reviewers for their valuable comments and suggestions which 
greatly improved the presentation of this paper.

%%%%%%%%%%%%%%%%%%
%%%%%%%%%%%%%%%%%%
%rrrrrrr
\newpage
%%%%%%%%%%%%%%%%%%
\def\cprime{$'$} \def\cprime{$'$} \def\cprime{$'$}
\providecommand{\bysame}{\leavevmode\hbox to3em{\hrulefill}\thinspace}
\providecommand{\MR}{\relax\ifhmode\unskip\space\fi MR }
% \MRhref is called by the amsart/book/proc definition of \MR.
\providecommand{\MRhref}[2]{%
  \href{http://www.ams.org/mathscinet-getitem?mr=#1}{#2}
}
\providecommand{\href}[2]{#2}

\end{document}